\input ./style/arxiv-general.cfg
\documentclass[aos,MSNbibl,nameyear,seceqn,dvips]{arximspdf}
\makeatletter
   \@ifpackageloaded{graphicx}{}{\usepackage{graphicx}}
\makeatother
\usepackage{multirow,dcolumn}

\doi{10.1214/15-AOS1320}
\volume{43}
\issue{4}
\pubyear{2015}
\firstpage{1682}
\lastpage{1715}
\docsubty{FLA}

\makeatletter
\newcolumntype{d}[1]{D{.}{.}{#1}}

\def\afrac#1#2{#1/(#2)}

\newcommand{\rrVert}{\Vert}
\newcommand{\llVert}{\Vert}
\renewcommand{\mid}{|}
\newtheorem{theorem}{Theorem}
\newproclaim{assumption}{Assumption}
\newproclaim{remark}{Remark}
\newtheorem{corollary}{Corollary}
\makeatother

\begin{document}
\begin{frontmatter}

\title{Semiparametric GEE analysis in partially linear single-index
models for longitudinal data}
\runtitle{Partially linear single-index models}

\begin{aug}
\author[A]{\fnms{Jia}~\snm{Chen}\thanksref{M1}\ead[label=e1]{jia.chen@york.ac.uk}},
\author[B]{\fnms{Degui}~\snm{Li}\thanksref{M1}\ead[label=e2]{degui.li@york.ac.uk}},
\author[C]{\fnms{Hua}~\snm{Liang}\thanksref{M2,T1}\ead[label=e3]{hliang@gwu.edu}}
\and
\author[D]{\fnms{Suojin}~\snm{Wang}\corref{}\thanksref{M3,T2}\ead[label=e4]{sjwang@stat.tamu.edu}}
\runauthor{Chen, Li, Liang and Wang}
\affiliation{University of York\thanksmark{M1}, George Washington
University\thanksmark{M2} and Texas A\&M University\thanksmark{M3}}
\address[A]{J. Chen\\
Department of Economics and Related Studies\\
University of York\\
Heslington West Campus\\
York, YO10 5DD\\
United Kingdom\\
\printead{e1}}
\address[B]{D. Li\\
Department of Mathematics\hspace*{5pt}\\
University of York\\
Heslington West Campus\\
York, YO10 5DD\\
United Kingdom\\
\printead{e2}}
\address[C]{H. Liang\\
Department of Statistics\\
George Washington University\\
Washington, District of Columbia 20052\\
USA\\
\printead{e3}}
\address[D]{S. Wang\\
Department of Statistics\\
Texas A\&M University\\
College Station, Texas 77843\\
USA\\
\printead{e4}}
\end{aug}
\thankstext{T1}{Supported in part by NSF Grants DMS-14-40121 and DMS-14-18042 and by Award Number 11228103, made
by National Natural Science Foundation of China.}
\thankstext{T2}{Supported in part by Award Number KUS-CI-016-04, made by King
Abdullah University of Science and Technology (KAUST).}

%
\received{\smonth{5} \syear{2014}}
%
\revised{\smonth{2} \syear{2015}}

%
\begin{abstract}
In this article, we study a partially linear single-index model for
longitudinal data under
a general framework which includes both the sparse and dense
longitudinal data cases. A semiparametric estimation method based on a
combination of the local linear smoothing and generalized estimation
equations (GEE) is introduced to estimate the two parameter vectors as
well as the unknown link function. Under some mild conditions, we
derive the asymptotic properties of the proposed parametric and
nonparametric estimators in different scenarios, from which we find
that the convergence rates and asymptotic variances of the proposed
estimators for sparse longitudinal data would be substantially
different from those for dense longitudinal data. We also discuss the
estimation of the covariance (or weight) matrices involved in the
semiparametric GEE method. Furthermore, we provide some numerical
studies including Monte Carlo simulation and an empirical application
to illustrate our methodology and theory.
\end{abstract}

%
\end{frontmatter}


\section{Introduction}\label{sec1}

Consider a semiparametric partially linear single-index model defined by
%
\begin{equation}
\label{eq11} Y(t)={\mathbf{Z}}^\top(t){\bolds\beta}+\eta\bigl({
\mathbf{X}}^{\top}(t){\bolds\theta} \bigr)+e(t),\qquad t\in{\mathcal{T}},
\end{equation}
where ${\mathcal{T}}$ is a bounded time interval, ${\bolds\beta
}$ and ${\bolds\theta}$ are two unknown vectors of parameters
with dimensions $d$ and $p$, respectively, $\eta(\cdot)$ is an
unknown link function, $Y(t)$ is a scalar stochastic process, ${\mathbf
{Z}}(t)$ and ${\mathbf{X}}(t)$ are covariates with dimensions $d$ and~$p$, respectively, and $e(t)$ is the random error process. For the case
of independent and identically distributed (i.i.d.) or weakly dependent
time series data, there has been extensive literature on statistical
inference of model (\ref{eq11}) since its introduction by \citet{CFGW97}. Several different approaches have been proposed to
estimate the unknown parameters and link function involved; see, for
example, \citet{XTL99}, \citet{YR02}, \citet{XH06}, \citet{WXZC10} and \citet{MZ13}. A recent paper by
\citet{LLLT10} further developed semiparametric techniques for the
variable selection and model specification testing issues in the
context of model (\ref{eq11}).

In this paper, we are interested in studying partially linear
single-index model (\ref{eq11}) in the context of longitudinal data
which arise frequently in many fields of research, such as biology,
climatology, economics and epidemiology, and thus have attracted
considerable attention in the literature in recent years. Various
parametric models and methods have been studied in depth for
longitudinal data; see \citet{DHLZ02} and the references therein.
However, the parametric models may be misspecified in practice, and the
misspecification may lead to inconsistent estimates and incorrect
conclusions being drawn. Hence, to circumvent this issue, in recent
years, there has been a large literature on how to relax the parametric
assumptions on longitudinal data models and many nonparametric, and
semiparametric models have thus been investigated; see, for example,
\citet{LY01}, \citet{HeZhuFun02}, \citet{FL04}, \citet{WCL05},
\citet{LC06}, \citet{WZ06}, \citet{LH10}, \citet{JW11} and \citet{YL13}.

Suppose that we have a random sample with $n$ subjects from model (\ref
{eq11}). For the $i$th subject, $i=1,\dots,n$, the response variable
$Y_i(t)$ and the covariates $ \{{\mathbf{Z}}_i(t), {\mathbf
{X}}_i(t) \}$ are collected at random time points $t_{ij}$,
$j=1,\dots,m_i$, which are distributed in a bounded time interval
${\mathcal{T}}$ according to the probability density function
$f_T(t)$. Here $m_i$ is the total number of observations for the $i$th
subject. To accommodate such longitudinal data, model (\ref{eq11}) is
written in the following framework:
%
\begin{equation}
\label{eq12} Y_i(t_{ij})={\mathbf{Z}}_i^\top(t_{ij}){
\bolds\beta}+\eta\bigl({\mathbf{X}}_i^{\top}(t_{ij}){
\bolds\theta} \bigr)+e_i(t_{ij})
\end{equation}
for $i=1,\dots,n$ and $j=1,\dots,m_i$. When $m_i$ varies across the
subjects, the longitudinal data set under investigation is unbalanced.
Several nonparametric and semiparametric models can be viewed as
special cases of model (\ref{eq12}). For instance, when ${\bolds
\beta}={\mathbf0}$, model (\ref{eq12}) reduces to the single-index
longitudinal data model [\citet{JW11}, \citet{CGL13a}];
when $p=1$ and ${\bolds\theta}=1$, model (\ref{eq12}) reduces
to the partially linear longitudinal data model [\citet{FL04}]. To
avoid confusion, we let ${\bolds\beta}_0$ and ${\bolds
\theta}_0$ be the true values of the two parameter vectors. For
identifiability reasons, ${\bolds\theta}_0$ is assumed to be a
unit vector with the first nonzero element being positive.
Furthermore, we allow that there exists certain within-subject
correlation structure for $e_i(t_{ij})$, which makes the model
assumption more realistic but the development of estimation methodology
more challenging.

To estimate the parameters ${\bolds\beta}_0$, ${\bolds
\theta}_0$ as well as the link function $\eta(\cdot)$ in model~(\ref
{eq12}), we first apply the local linear approximation to the unknown
link function, and then introduce a profile weighted least squares
approach to estimate the two parameter vectors based on the technique
of generalized estimation equations (GEE). Under some mild conditions,
we derive the asymptotic properties of the developed parametric and
nonparametric estimators in different scenarios. Our framework is
flexible in that $m_i$ can either be bounded or tend to infinity. Thus
both the dense and sparse longitudinal data cases can be included.
Dense longitudinal data means that there exists a sequence of positive
numbers $M_n$ such that $\min_{i}m_i\geq M_n$, and $M_n\rightarrow
\infty$ as $n\rightarrow\infty$ [see, e.g., \citet{HMW06} and
\citet{ZC07}], whereas sparse longitudinal data means that
there exists a positive constant $M_*$ such that $\max_{i}m_i\leq
M_*$; see, for example, \citet{YMW}, \citet{WQC10}. We show
that the convergence rates and asymptotic variances of our
semiparametric estimators in the sparse case are substantially
different from those in the dense case. Furthermore, we show that the
proposed semiparametric GEE (SGEE)-based estimators are asymptotically
more efficient than the profile unweighted least squares (PULS)
estimators, when the weights in the SGEE method are chosen as the
inverse of the covariance matrix of the errors.
We also introduce a semiparametric approach to estimate the
covariance matrices (or weights) involved in the SGEE method, which is
based on a variance--correlation decomposition and consists of two
steps: first, estimate the conditional variance function using a robust
nonparametric method that accommodates heavy-tailed errors, and second,
estimate the parameters in the correlation matrix. A simulation study
and a real data analysis are provided to illustrate our methodology and theory.

The rest of the paper is organized as follows. In Section~\ref{sec2}, we
introduce the SGEE methodology for estimating ${\bolds\beta}_0$,
${\bolds\theta}_0$ and $\eta(\cdot)$. Section~\ref{sec3} establishes
the large sample theory for the proposed parametric and nonparametric
estimators and gives some related discussions. Section~\ref{sec4} discusses how
to determine the weight matrices in the estimation equations. Section~\ref{sec5}
gives some numerical examples to investigate the finite sample
performance of the proposed approach. Section~\ref{sec6} concludes the paper.
Technical assumptions are given in Appendix~\ref{appA}. The proofs of the main
results are given in Appendix~\ref{appB}. Some auxiliary lemmas and their proofs
are provided in the supplementary material [\citet{CLLWsupp}].


\section{Estimation methodology}
\label{sec2}

Various semiparametric estimation\break approaches have been proposed to
estimate model (\ref{eq11}) in the case of i.i.d.~observations (or
weakly dependent time series data). See, for example, \citet{CFGW97} and \citet{LLLT10} for the profile
likelihood method, \citet{YR02} and \citet{WXZC10} for the ``remove-one-component'' technique using penalized spline
and local linear smoothing, respectively, and \citet{XH06}
for the minimum average variance estimation approach. However, there
is limited literature on partially linear single-index models for
longitudinal data because of the more complicated structures
involved. Recently, \citet{CGL13b} studied a partially
linear single-index longitudinal data model with individual effects.
To remove the individual effects and derive consistent
semiparametric estimators, they had to limit their discussions to
the dense and balanced longitudinal data case.
\citet{MLT} considered a partially linear single-index
longitudinal data
model by using polynomial splines to approximate the unknown link
function, but their discussion was limited to the sparse and
balanced longitudinal data case. In contrast, as mentioned in
Section~\ref{sec1}, our framework includes both the sparse and dense
longitudinal data cases. Meanwhile, observations are allowed to be
collected at irregular and subject specific time points. All this
provides much wider applicability of our framework. Furthermore, to
improve the efficiency of the semiparametric estimation, we develop
a new profile weighted least squares approach to estimate the
parameters ${\bolds\beta}_0$, ${\bolds\theta}_0$ as well
as the link function $\eta_0(\cdot)$.

To simplify the presentation, let
\begin{eqnarray*}
{\mathbf Y}_i&=& \bigl(Y_{i}(t_{i1}),
\dots,Y_{i}(t_{im_i}) \bigr)^\top,\qquad {
\mathbf{X}}_i= \bigl({\mathbf{X}}_{i}(t_{i1}),
\dots,{\mathbf{X}}_{i}(t_{im_i}) \bigr)^\top,
\\
{\mathbf{Z}}_i&=& \bigl({\mathbf{Z}}_{i}(t_{i1}),
\dots,{\mathbf{Z}}_{i}(t_{im_i}) \bigr)^\top,\qquad {
\mathbf e}_i= \bigl(e_{i}(t_{i1}),
\dots,e_{i}(t_{im_i}) \bigr)^\top,
\\
{\bolds\eta}({\mathbf{X}}_i,{\bolds\theta})&=& \bigl(\eta
\bigl({\mathbf{X}}_{i}^\top(t_{i1}){\bolds\theta} \bigr),\dots, \eta\bigl({\mathbf{X}}_{i}^\top(t_{im_i}){
\bolds\theta} \bigr) \bigr)^\top.
\end{eqnarray*}
With the above notation, model (\ref{eq12}) can then be re-written as
%
\begin{equation}
\label{eq21} {\mathbf Y}_i={\mathbf{Z}}_i{\bolds\beta}_0+{\bolds\eta}({\mathbf{X}}_i,{\bolds\theta}_0)+{\mathbf e}_i.
\end{equation}
We further let
${\mathbb Y}=({\mathbf Y}_1^\top,\dots,{\mathbf Y}_n^\top)^\top$,
${\mathbb{Z}}= ({\mathbf{Z}}_1^\top,\dots,{\mathbf{Z}}_n^\top
)^\top$, ${\mathbb{E}}=({\mathbf e}_1^\top,\dots,{\mathbf
e}_n^\top)^\top$,
${\bolds\eta}({\mathbb{X}},{\bolds\theta})=
({\bolds\eta}^\top({\mathbf{X}}_1, {\bolds\theta
} ),\dots,
{\bolds\eta}^\top({\mathbf{X}}_{n},{\bolds\theta})
)^\top$.
Then model (\ref{eq21}) is equivalent to
%
\begin{equation}
\label{eq22} {\mathbb Y}={\mathbb Z} {\bolds\beta}_0+{
\bolds\eta}({\mathbb{X}},{\bolds\theta}_0)+{\mathbb E}.
\end{equation}

Our estimation procedure is based on the profile likelihood method,
which is commonly used in semiparametric estimation; see, for example, \citet{CFGW97}, \citet{FH05} and
\citet{FHL07}. Let
$Y_{ij}=Y_{i}(t_{ij})$, ${\mathbf{Z}}_{ij}={\mathbf{Z}}_{i}(t_{ij})$
and ${\mathbf{X}}_{ij}={\mathbf{X}}_{i}(t_{ij})$. For given
${\bolds\beta}$ and ${\bolds\theta}$, we can estimate
$\eta(\cdot)$ and its derivative $\dot{\eta}(\cdot)$ at point $u$
by minimizing the following loss function:
%
\begin{eqnarray}\label{eq23}
&& L_n (a,b\mid{\bolds\beta},{\bolds\theta} )
\nonumber\\[-8pt]\\[-8pt]\nonumber
&&\qquad =
\sum_{i=1}^n \Biggl\{\frac{w_i}{h}\sum
_{j=1}^{m_i} \bigl[Y_{ij}-{\mathbf
{Z}}_{ij}^\top{\bolds\beta}-a-b\bigl({
\mathbf{X}}_{ij}^\top{\bolds\theta}-u\bigr)
\bigr]^2 K \biggl(\frac{{\mathbf{X}}_{ij}^\top{\bolds\theta}-u}{h}
\biggr) \Biggr\},
\end{eqnarray}
where $K(\cdot)$ is a kernel function, $h$ is a bandwidth and $w_i$,
$i=1,\dots,n$, are some weights. It is well known that the local
linear smoothing has advantages over the Nadaraya--Watson kernel
method, such as higher asymptotic efficiency, design adaption and
automatic boundary correction [\citet{FG96}]. Following the
existing literature such as \citet{WZ06}, the weights $w_i$ can
be specified by two schemes: $w_i=1/T_n$ (type 1) and $w_i=1/(nm_i)$
(type 2), where $T_n=\sum_{i=1}^nm_i$. The type 1 weight scheme
corresponds to an equal weight for each observation, while the type 2
scheme corresponds to an equal weight within each subject. As discussed
in \citet{HWZ02} and \citet{WZ06}, the type 2 scheme may
be appropriate if the number of observations varies across subjects. As
the longitudinal data under investigation in this paper are allowed to
be unbalanced, we use $w_i=1/(nm_i)$, which was also used by \citet{LH10} and \citet{KZ12}. We denote
%
\begin{equation}
\label{eq24} \bigl(\widehat{\eta}(u\mid{\bolds\beta},{\bolds\theta
}), \widehat{\dot{\eta}}(u\mid{\bolds\beta},{\bolds\theta})
\bigr)^\top=\arg\min_{a,b}L_n (a,b\mid{
\bolds\beta},{\bolds\theta} ).
\end{equation}
By some elementary calculations [see, e.g., \citet{FG96}], we have
%
\begin{equation}
\label{eq25} \widehat{\eta}(u\mid{\bolds\beta}, {\bolds\theta})=\sum
_{i=1}^n{\mathbf s}_i(u\mid{
\bolds\theta}) ({\mathbf Y}_i-{\mathbf{Z}}_i{
\bolds\beta} )
\end{equation}
for given ${\bolds\beta}$ and ${\bolds\theta}$, where
\begin{eqnarray}
\qquad {\mathbf s}_i(u\mid{\bolds\theta})&=&(1,0) \Biggl[\sum
_{i=1}^n\overline{\mathbf{X}}_i^\top(u
\mid{\bolds\theta}){\mathbf K}_{i} (u\mid{\bolds\theta})
\overline{\mathbf{X}}_i(u\mid{\bolds\theta})
\Biggr]^{-1}\overline{\mathbf{X}}_i^\top(u\mid{
\bolds\theta}){\mathbf K}_{i} (u\mid{\bolds\theta}),
\nonumber
\\
\overline{\mathbf{X}}_i(u\mid{\bolds\theta})&=& \bigl(
\overline{\mathbf{X}}_{i1}(u\mid{\bolds\theta}),\dots,\overline{
\mathbf{X}}_{im_i}(u\mid{\bolds\theta}) \bigr)^\top,
\nonumber\\[-8pt]\\[-8pt]\nonumber
\overline{\mathbf{X}}_{ij}(u\mid{\bolds\theta})&=& \bigl(1,{
\mathbf{X}}_{ij}^\top{\bolds\theta}-u
\bigr)^\top,
\\
{\mathbf K}_{i}(u\mid{\bolds\theta})&=&\operatorname{diag}
\biggl(w_iK \biggl(\frac{{\mathbf{X}}_{i1}^\top{\bolds\theta
}-u}{h} \biggr),\dots,w_iK
\biggl(\frac{{\mathbf{X}}_{im_i}^\top
{\bolds\theta}-u}{h} \biggr) \biggr).
\nonumber
\end{eqnarray}

Based on the profile least squares approach with the first-stage local
linear smoothing, we can construct estimators of the parameters
${\bolds\beta}_0$ and ${\bolds\theta}_0$. We start with
the PULS method which ignores the possible within-subject correlation
structure. Define the PULS loss function by
%
\begin{eqnarray}\label{eq26}
Q_{n0}({\bolds\beta},{\bolds\theta})&=&\sum
_{i=1}^n \bigl[{\mathbf Y}_i-{
\mathbf{Z}}_i{\bolds\beta}-\widehat{\bolds\eta}({
\mathbf{X}}_i\mid{\bolds\beta}, {\bolds\theta})
\bigr]^\top\bigl[{\mathbf Y}_i-{\mathbf{Z}}_i{
\bolds\beta}-\widehat{\bolds\eta}({\mathbf{X}}_i\mid{
\bolds\beta}, {\bolds\theta}) \bigr]
\nonumber\\[-8pt]\\[-8pt]\nonumber
&=& \bigl[{\mathbb Y}-{\mathbb{Z}} {\bolds\beta}-\widehat
{\bolds\eta}({\mathbb{X}}\mid{\bolds\beta}, {\bolds\theta})
\bigr]^\top\bigl[{\mathbb Y}-{\mathbb{Z}} {\bolds\beta}-\widehat{
\bolds\eta}({\mathbb{X}}\mid{\bolds\beta}, {\bolds\theta})
\bigr],
\end{eqnarray}
where, for given ${\bolds\beta}$ and ${\bolds\theta}$,
$\widehat{\bolds\eta}({\mathbf{X}}_i|{\bolds\beta},
{\bolds\theta})$ and $\widehat{\bolds\eta}({\mathbb
{X}}|{\bolds\beta}, {\bolds\theta})$ are the local linear
estimators of the vectors ${\bolds\eta}({\mathbf
{X}}_i,{\bolds\theta})$ and ${\bolds\eta}({\mathbb
{X}},{\bolds\theta})$, respectively; that is, each element of
$\widehat{\bolds\eta}({\mathbf{X}}_i|{\bolds\beta},
{\bolds\theta})$ and $\widehat{\bolds\eta}({\mathbb
{X}}|{\bolds\beta}, {\bolds\theta})$ is defined as in
(\ref{eq25}). The PULS estimators of ${\bolds\beta}_0$ and
${\bolds\theta}_0$ are obtained by minimizing the loss function
$Q_{n0}({\bolds\beta},{\bolds\theta})$ with\vspace*{1pt} respect to
${\bolds\beta}$ and ${\bolds\theta}$ and normalizing the
minimizer ${\bolds\theta}$. We denote the resulting estimators
by $\widetilde{\bolds\beta}$ and~$\widetilde{\bolds\theta}$, respectively.

Although it is easy to verify that both $\widetilde{\bolds\beta
}$ and $\widetilde{\bolds\theta}$ are consistent, they are not
efficient as the within-subject correlation structure is not taken into
account. Hence, to improve the efficiency of the parametric estimators,
we next introduce a GEE-based method to estimate the parameters
${\bolds\beta}_0$ and ${\bolds\theta}_0$. Existing
literature on \mbox{GEE-}based method in longitudinal data analysis includes
\citet{LZ86}, \citet{XY03} and \citet{W11}. Let
${\mathbb{W}}=\operatorname{diag}\{{\mathbf{W}}_1,\dots,{\mathbf{W}}_n\}$,
where ${\mathbf{W}}_i={\mathbf{R}}_i^{-1}$ and ${\mathbf{R}}_i$ is
an $m_i\times m_i$ working covariance matrix whose estimation will be
discussed in Section~\ref{sec4}. Define
\begin{eqnarray}
{\bolds\rho}_{\mathbf{Z}}({\mathbf{X}}_{i},{\bolds\theta}) &=& \bigl(\rho_{\mathbf{Z}}\bigl({\mathbf{X}}_{i1}^\top
{\bolds\theta}\mid{\bolds\theta}\bigr),\dots, {\rho}_{\mathbf{Z}}
\bigl({\mathbf{X}}_{im_i}^\top{\bolds\theta}\mid{
\bolds\theta}\bigr) \bigr)^\top,\qquad \rho_{\mathbf
{Z}}(u\mid{
\bolds\theta})=\mathrm{E} \bigl[{\mathbf{Z}}_{ij}\mid{\mathbf
{X}}_{ij}^\top{\bolds\theta}=u \bigr],
\nonumber
\\
{\bolds\rho}_{\mathbf{X}}({\mathbf{X}}_{i},{\bolds\theta})&=& \bigl(\rho_{\mathbf{X}}\bigl({\mathbf{X}}_{i1}^\top
{\bolds\theta}\mid{\bolds\theta}\bigr),\dots, \rho_{\mathbf{X}}
\bigl({\mathbf{X}}_{im_i}^\top{\bolds\theta}\mid{
\bolds\theta}\bigr) \bigr)^\top,\qquad \rho_{\mathbf
{X}}(u\mid{
\bolds\theta})=\mathrm{E} \bigl[{\mathbf{X}}_{ij}\mid{\mathbf
{X}}_{ij}^\top{\bolds\theta}=u \bigr],
\nonumber
\\
{\bolds\Lambda}_i({\bolds\theta}) &=& \bigl({\mathbf
{Z}}_i-{\bolds\rho}_{\mathbf{Z}}({\mathbf{X}}_i,{
\bolds\theta}), \bigl[\dot{\bolds\eta}({\mathbf{X}}_i,{
\bolds\theta})\otimes{\mathbf1}_{p}^\top\bigr]\odot
\bigl[{\mathbf{X}}_i-{\bolds\rho}_{\mathbf X}({
\mathbf{X}}_i,{\bolds\theta}) \bigr] \bigr),
\nonumber
\end{eqnarray}
where $\dot{\bolds\eta}({\mathbf{X}}_i,{\bolds\theta})$
is a column\vspace*{1pt} vector with its elements being the derivatives of $\eta
(\cdot)$ at points ${\mathbf{X}}_{ij}^\top{\bolds\theta}$,
$j=1, \dots, m_i$, ${\mathbf1}_{p}$ is a $p$-dimensional\vspace*{1pt} vector
of ones, $\otimes$ is the Kronecker product and $\odot$ denotes the
componentwise product. The construction of the parametric estimators is
based on solving the following equation with respect to ${\bolds\beta}$ and ${\bolds\theta}$:
%
\begin{equation}
\label{eq27} \sum_{i=1}^n\widehat{
\bolds\Lambda}_i^\top({\bolds\theta}) {
\mathbf{W}}_i \bigl[{\mathbf Y}_i-{\mathbf{Z}}_i{
\bolds\beta}-\widehat{\bolds\eta}({\mathbf{X}}_i\mid{
\bolds\beta}, {\bolds\theta}) \bigr]={\mathbf{0}},
\end{equation}
where $\widehat{\bolds\Lambda}_i({\bolds\theta})$ is an
estimator of ${\bolds\Lambda}_i({\bolds\theta})$ with
${\bolds\rho}_{\mathbf{Z}}({\mathbf{X}}_i,{\bolds\theta
})$, ${\bolds\rho}_{\mathbf{X}}({\mathbf{X}}_i,{\bolds\theta})$ and $\dot{\bolds\eta}({\mathbf{X}}_i,{\bolds\theta})$ replaced by their corresponding local linear estimated
values. Let $\widehat{\bolds\beta}$ and $\widehat{\bolds\theta}_1$ be the solutions to the estimation equations in (\ref
{eq27}), and let the SGEE-based estimator of ${\bolds\theta}_0$
be defined as $\widehat{\bolds\theta}=\widehat{\bolds\theta}_1/\llVert \widehat{\bolds\theta}_1 \rrVert $, where
$\llVert
\cdot\rrVert $ is the Euclidean norm. Note that the solutions to the
equations in (\ref{eq27}) generally do not have a closed form. In the
numerical studies, we use the trust-region dogleg algorithm within the
Matlab command ``fsolve'' to obtain the solutions to (\ref{eq27}).
Corollary~\ref{co1} below shows that the SGEE-based estimators $\widehat
{\bolds\beta}$ and $\widehat{\bolds\theta}$ are
generally asymptotically more efficient than the PULS estimators
$\widetilde{\bolds\beta}$ and $\widetilde{\bolds\theta
}$, when the weights are chosen appropriately.

Replacing ${\bolds\beta}$ and ${\bolds\theta}$ in
$\widehat{\eta}(\cdot)$ by $\widehat{\bolds\beta}$ and
$\widehat{\bolds\theta}$, respectively, we obtain the local
linear estimator of the link function $\eta(\cdot)$ at $u$ as
%
\begin{equation}
\label{eq28} \widehat{\eta}(u)=\widehat{\eta}(u\mid\widehat{\bolds\beta
},\widehat{\bolds\theta})=\sum_{i=1}^n{
\mathbf s}_i(u \mid\widehat{\bolds\theta}) ({\mathbf
Y}_i-{\mathbf{Z}}_i\widehat{\bolds\beta} ).
\end{equation}

In Section~\ref{sec3} below, we will give the large sample properties of the
estimators proposed above, and in Section~\ref{sec4}, we will discuss how to
choose the working covariance matrix ${\mathbf{R}}_i$.


\section{Theoretical properties}\label{sec3}

Before establishing the large sample theory for the proposed parametric
and nonparametric estimators, we introduce some notation. Let ${\mathbf
B}_0$ be a $p\times(p-1)$ matrix such that ${\mathbf M}=
({\bolds\theta}_0, {\mathbf B}_0)$ is a $p\times p$ orthogonal
matrix, and define
\[
{\mathbf I}({\mathbf B}_0)=\pmatrix{ {\mathbf I}_d & {
\mathbf O}_{d\times(p-1)}
\vspace*{3pt}\cr
{\mathbf O}_{p\times d} & {\mathbf B}_0},
\]
where ${\mathbf I}_k$ is a $k\times k$ identity matrix and ${\mathbf
O}_{k\times l}$ is a $k\times l$ null matrix. Let ${\bolds\Lambda
}_i={\bolds\Lambda}_i({\bolds\theta}_0)$, and assume that
there exist two positive semi-definite matrices ${\bolds\Omega
}_0$ and ${\bolds\Omega}_1$ as well as a sequence of numbers
$\omega_n$ such that $\omega_n\rightarrow\infty$,
%
\begin{eqnarray}
\frac{1}{\omega_n}\sum_{i=1}^n{
\bolds\Lambda}_i^\top{\mathbf{W}}_i{
\bolds\Lambda}_i & \stackrel{P}\rightarrow&{\bolds\Omega}_0,\label{eq31}
\\
\frac{1}{\omega_n}\sum_{i=1}^n\mathrm{E}
\bigl[{\bolds\Lambda}_i^\top{\mathbf{W}}_i{
\mathbf e}_i{\mathbf e}_i^\top{\mathbf
{W}}_i{\bolds\Lambda}_i \bigr]&\rightarrow&{
\bolds\Omega}_1,\label{eq32}
\\
\max_{1\leq i\leq n}\mathrm{E} \bigl[{\bolds\Lambda}_i^\top
{\mathbf{W}}_i{\mathbf e}_i{\mathbf
e}_i^\top{\mathbf{W}}_i{\bolds\Lambda}_i \bigr]&=&o(\omega_n),\label{eq321}
\end{eqnarray}
as $n\rightarrow\infty$, and ${\mathbf I}^\top({\mathbf
B}_0){\bolds\Omega}_0 {\mathbf I}({\mathbf B}_0)$ is positive
definite. Conditions (\ref{eq32}) and (\ref{eq321}) ensure that
the Lindeberg--Feller condition can be satisfied, and thus the
classical central limit theorem for independent sequence [\citet{P95}] is applicable. It is not difficult to verify the assumption in
(\ref{eq321}) for the dense and sparse longitudinal data. In
particular, (\ref{eq321}) excludes the case where the term
${\bolds\Lambda}_i^\top{\mathbf{W}}_i{\mathbf e}_i$ from one
or a few subjects dominates those from the others. For the latter case,
it may be possible to derive the consistency of the proposed parametric
estimation, but the proof of the asymptotic normality would be
difficult. Let ${\bolds\Omega}_0^{+}$ be the Moore--Penrose
inverse matrix of ${\bolds\Omega}_0$, which is defined as
${\bolds\Omega}_0^{+}={\mathbf I}({\mathbf B}_0) [{\mathbf
I}^\top({\mathbf B}_0){\bolds\Omega}_0 {\mathbf I}({\mathbf
B}_0) ]^{-1}{\mathbf I}^\top({\mathbf B}_0)$. We next give the
asymptotic distribution theory for the SGEE-based estimators $\widehat
{\bolds\beta}$ and $\widehat{\bolds\theta}$.

\begin{theorem}\label{teo1}
Suppose that Assumptions~\ref{as1}--\ref{as5} in Appendix~\ref{appA} and (\ref{eq31})--(\ref
{eq321}) are satisfied. Then we have
%
\begin{equation}
\label{eq33} \omega_n^{1/2} \pmatrix{ \widehat{\bolds\beta}-{\bolds\beta}_0
\vspace*{3pt}\cr
\widehat{\bolds\theta}-{
\bolds\theta}_0} \stackrel{d}\longrightarrow \mathrm{N} \bigl({
\mathbf0}, {\bolds\Omega}_0^{+}{\bolds\Omega}_1{\bolds\Omega}_0^{+} \bigr)
\end{equation}
as $n\rightarrow\infty$.
\end{theorem}

%
\begin{remark}\label{re1}
{Theorem~\ref{teo1} establishes the asymptotically normal distribution
theory for $\widehat{\bolds\beta}$ and $\widehat{\bolds\theta}$ with convergence rate $\omega_n^{1/2}$. This $\omega_n$ is
linked to $h$ through $n$ in a certain way. Specifically, the condition
$\omega_nh^6\rightarrow0$ in Assumption~\ref{as5} needs to be satisfied to
ensure that the bias term of the parametric estimation is
asymptotically negligible. The specific forms of $\omega_n$,
${\bolds\Omega}_0$ and ${\bolds\Omega}_1$ can be derived
for some particular cases, for instance, when longitudinal data are
balanced, that is, $m_i\equiv m$, $\omega_n=nm$. Furthermore, assume
that the covariates and the error are i.i.d. with $\mathrm{E}
[e_i^2(t_{ij}) ]\equiv\sigma_e^2$, $e_i(t_{ij})$ is independent
of the covariates and ${\mathbf{W}}_i$, $i=1,\dots,n$, are $m\times
m$ identity matrices. Then we can show that
\[
{\bolds\Omega}_0=\pmatrix{ \Omega_0(1)&
\Omega_0(2)
\vspace*{3pt}\cr
\Omega_0^\top(2)&
\Omega_0(3)}\quad\mbox{and}\quad{\bolds\Omega}_1=
\sigma_e^2\pmatrix{ \Omega_0(1)&
\Omega_0(2)
\vspace*{3pt}\cr
\Omega_0^\top(2)&
\Omega_0(3)},
\]
where
\begin{eqnarray}
\Omega_0(1)&=&\mathrm{E} \bigl\{ \bigl[{\mathbf{Z}}(t)-{\bolds\rho
}_{\mathbf{Z}}\bigl({\mathbf{X}}^\top(t){\bolds\theta
}_0\mid{\bolds\theta}_0\bigr) \bigr] \bigl[{
\mathbf{Z}}(t)-{\bolds\rho}_{\mathbf{Z}}\bigl({\mathbf{X}}^\top(t){
\bolds\theta}_0\mid{\bolds\theta}_0\bigr)
\bigr]^\top\bigr\},
\nonumber
\\
\Omega_0(2)&=&\mathrm{E} \bigl\{\dot{\eta}\bigl({\mathbf{X}}^\top
(t){\bolds\theta}_0\bigr) \bigl[{\mathbf{Z}}(t)-
\rho_{\mathbf{Z}}\bigl({\mathbf{X}}^\top(t){\bolds\theta}_0\mid{\bolds\theta}_0\bigr) \bigr] \bigl[{
\mathbf{X}}(t)-\rho_{\mathbf{X}}\bigl({\mathbf{X}}^\top(t){
\bolds\theta}_0\mid{\bolds\theta}_0\bigr)
\bigr]^\top\bigr\},
\nonumber
\\
\Omega_0(3)&=&\mathrm{E} \bigl\{ \bigl[\dot{\eta}\bigl({
\mathbf{X}}^\top(t){\bolds\theta}_0\bigr)
\bigr]^2 \bigl[{\mathbf{X}}(t)-\rho_{\mathbf{X}}\bigl({
\mathbf{X}}^\top(t){\bolds\theta}_0\mid{\bolds\theta}_0\bigr) \bigr] \bigl[{\mathbf{X}}(t)-\rho_{\mathbf{X}}
\bigl({\mathbf{X}}^\top(t){\bolds\theta}_0\mid{
\bolds\theta}_0\bigr) \bigr]^\top\bigr\}.
\nonumber
\end{eqnarray}
Hence ${\bolds\Omega}_0^{+}{\bolds\Omega}_1{\bolds\Omega}_0^{+}$ reduces to $\sigma_e^2{\bolds\Omega}_0^{+}$.

In Theorem~\ref{teo1} above, we only require $n\rightarrow\infty$. As
mentioned in Section~\ref{sec1}, both the sparse and dense
longitudinal data cases can be included in a unified framework. For the
sparse longitudinal data case when $m_i$ is bounded by a certain
positive constant, we can take $\omega_n=n$ and prove that (\ref
{eq33}) holds. For the dense longitudinal data case where $\min_i
m_i\geq M_n$ with $M_n\rightarrow\infty$, under some regularity conditions
we may prove (\ref{eq33}) with $w_n=\sum_{i=1}^n m_i$.
As more observations are available in the dense longitudinal data case
and the order for the total number of the observations is higher than
$n$, the convergence rate for the parametric estimators is faster than
the well-known root-$n$ rate in the sparse longitudinal data case. }
\end{remark}

Using Theorem~\ref{teo1}, we can obtain the following corollary.

\begin{corollary}\label{co1}
Suppose that the weights ${\mathbf{W}}_i$ in (\ref{eq27}) are chosen
as the inverse of the conditional covariance matrix of ${\mathbf
{e}}_i$, and the
conditions of Theorem~\ref{teo1} are satisfied. Then the SGEE-based estimators
$\widehat{\bolds\beta}$ and $\widehat{\bolds\theta}$
are asymptotically more efficient than the PULS estimators $\widetilde
{\bolds\beta}$ and $\widetilde{\bolds\theta}$ defined in
Section~\ref{sec2}.
\end{corollary}

\begin{remark}\label{re2}
In\vspace*{1pt} the proof of the above corollary, we show that the asymptotic
covariance matrix of the PULS estimators $\widetilde{\bolds\beta
}$ and $\widetilde{\bolds\theta}$ (after appropriate
normalization) minus that of the SGEE-based estimators $\widehat
{\bolds\beta}$ and $\widehat{\bolds\theta}$ is positive
semi-definite, although the two estimation methods have the same
convergence rates. That is, under the conditions assumed in Theorem~\ref{teo1},
the limit matrix of $\omega_n [\operatorname{Var}(\widetilde{\bolds\beta}, \widetilde{\bolds\theta})-\operatorname{Var}(\widehat
{\bolds\beta}, \widehat{\bolds\theta}) ]$ is
positive semi-definite. For the case of independent observations, a~recent paper by \citet{LBY14} discussed the efficient bound for the
semiparametric estimation in single-index models. Following their idea,
we conjecture that modification of our estimation procedure may be
needed to obtain the efficient estimation in the partially linear
single-index longitudinal data models. We will study this issue in our
future research.
\end{remark}

To establish the asymptotic distribution theory for the nonparametric
estimator $\widehat{\eta}(u)$ under a unified framework, we assume
that there exist a sequence $\varphi_n(h)$ and a constant $0<\sigma
_*^2<\infty$ such that
%
\begin{equation}
\label{eq340} \varphi_n(h)=o(\omega_n),\qquad {
\varphi_n(h)}\max_{1\leq i\leq n}\mathrm{E} \bigl[{\mathbf
s}_i(u\mid{\bolds\theta}_0){\mathbf
e}_i{\mathbf e}_i^\top{\mathbf
s}_i^\top(u\mid{\bolds\theta}_0)
\bigr]=o(1)
\end{equation}
and
%
\begin{equation}
\label{eq34} {\varphi_n(h)}\sum_{i=1}^n\mathrm{E} \bigl[{\mathbf s}_i(u\mid{\bolds\theta}_0){
\mathbf e}_i{\mathbf e}_i^\top{\mathbf
s}_i^\top(u\mid{\bolds\theta}_0)
\bigr]\rightarrow\sigma_*^2.
\end{equation}
The first restriction in (\ref{eq340}) is imposed to ensure that the
parametric convergence rates are faster than the nonparametric
convergence rates, and the second restriction in (\ref{eq340}) and
the condition in (\ref{eq34}) are imposed for the derivation of the
asymptotic variance of the local linear estimator $\widehat{\eta}(u)$
and the satisfaction of the Lindeberg--Feller condition. The specific
forms of $\varphi_n(h)$ and $\sigma_*^2$ will be discussed in Remark~\ref{re3} below. Let $\mu_j=\int v^jK(v)\,dv$ for $j=0, 1, 2$ and $\ddot{\eta
}_0(\cdot)$ be the second-order derivative of $\eta_0(\cdot)$.

\begin{theorem}\label{theorem2}\label{teo2}
Suppose that the conditions of Theorem~\ref{teo1}, (\ref{eq340}) and (\ref
{eq34}) are satisfied. Then we have
%
\begin{equation}
\label{eq35} \varphi_n^{1/2}(h) \bigl[\widehat{\eta}(u)-
\eta_0(u)-b_\eta(u)h^2 \bigr]\stackrel{d}
\longrightarrow \mathrm{N}\bigl(0,\sigma_*^2\bigr),
\end{equation}
where $b_\eta(u)=\ddot{\eta}_0(u)\mu_2/2$.
\end{theorem}

\begin{remark}\label{re3}
{Theorem~\ref{teo2} provides the asymptotically normal distribution theory for
the nonparametric estimator $\widehat{\eta}(u)$ with a convergence
rate $O_P (\varphi_n^{-1/2}(h)+h^2 )$. The forms of $\varphi
_n(h)$ and $\sigma_*^2$ in Theorem~\ref{teo2} depend on the type of the
longitudinal data under study, that is, whether it is sparse or dense.
We can derive their specific forms for some particular cases. Consider,
for example, the case where $e_{i}(t_{ij})=v_i+\varepsilon_{ij}$, in\vspace*{1pt}
which $\varepsilon_{ij}$ are i.i.d. across both $i$ and $j$ with $\mathrm{E}[\varepsilon_{ij}]=0$ and $\mathrm{E}[\varepsilon_{ij}^2]=\sigma
^2_\varepsilon$, and $\{v_{i}\}$ is an i.i.d. sequence\vspace*{1pt} of random
variables with $\mathrm{E}[v_{i}]=0$ and $\mathrm{E}[v_{i}^2]=\sigma^2_v$
and is independent of $\{\varepsilon_{ij}\}$. In this case, we note that
\begin{eqnarray*}
&& \mathrm{E} \Biggl\{ \Biggl[\sum_{j=1}^{m_i}K
\biggl(\frac{{\mathbf
{X}}_{ij}^\top{\bolds\theta}_0-u}{h} \biggr)e_{ij} \Biggr]^2 \Biggr\}
\\
&&\qquad = \mathrm{E} \Biggl\{ \Biggl[\sum_{j=1}^{m_i}K
\biggl(\frac{{\mathbf
{X}}_{ij}^\top{\bolds\theta}_0-u}{h} \biggr) (v_i+\varepsilon_{ij})
\Biggr]^2 \Biggr\}
\\
&&\qquad = \sum_{j=1}^{m_i}\mathrm{E}
\biggl[K^2 \biggl(\frac{{\mathbf{X}}_{ij}^\top
{\bolds\theta}_0-u}{h} \biggr) (v_i+
\varepsilon_{ij})^2 \biggr]
\\
&&\quad\qquad{}  +\sum
_{j_1\neq j_2}\mathrm{E} \biggl[K \biggl(\frac{{\mathbf{X}}_{ij_1}^\top
{\bolds\theta}_0-u}{h} \biggr)
K \biggl(\frac{{\mathbf{X}}_{ij_2}^\top{\bolds\theta}_0-u}{h} \biggr) (v_i+
\varepsilon_{ij_1}) (v_i+\varepsilon_{ij_2}) \biggr]
\\
&&\qquad \sim m_ih\nu_0f_{{\bolds\theta}_0}(u) \bigl(
\sigma_v^2+\sigma_\varepsilon^2
\bigr)+m_i(m_i-1)h^2\mu_0^2f_{{\bolds\theta
}_0}^2(u)
\sigma_v^2,
\end{eqnarray*}
where $\nu_0=\int K^2(v)\,dv$ and $f_{{\bolds\theta}_0}(\cdot)$
is the probability density function of ${\mathbf X}_{ij}^\top{\bolds\theta}_0$.

For the sparse longitudinal data case,
$m_i(m_i-1)h^2\mu_0^2f_{{\bolds\theta}_0}^2(u) \sigma_v^2$ is
dominated by $m_ih\nu_0f_{{\bolds\theta}_0}(u)(\sigma
_v^2+\sigma_\varepsilon^2)$, as $m_i$ is bounded and
$h\rightarrow0$. Then, by Lemma 1 in the supplementary document [\citet{CLLWsupp}] and some
elementary calculations, we can prove that
%
\begin{eqnarray}
\sum_{i=1}^n\mathrm{E} \bigl[{\mathbf
s}_i(u\mid{\bolds\theta}_0){\mathbf
e}_i{\mathbf e}_i^\top{\mathbf
s}_i^\top(u\mid{\bolds\theta}_0)
\bigr]&\sim&\frac{1}{(nh)^2}\sum_{i=1}^n
\frac{m_ih\nu
_0(\sigma_v^2+\sigma_\varepsilon
^2)}{m_i^2f_{{\bolds\theta}_0}(u)}
\nonumber\\[-8pt]\\[-8pt]\nonumber
&\sim&\frac{\nu_0(\sigma_v^2+\sigma
_\varepsilon^2)}{n^2hf_{{\bolds\theta}_0}(u)}\sum_{i=1}^n
\frac{1}{m_i}.\label{eq36}
\end{eqnarray}
Hence, in this case, we can take $\varphi_n(h)=(n^2h) (\sum
_{i=1}^n1/m_i )^{-1}$ which has the same order as $nh$, and $\sigma
_*^2=\nu_0(\sigma_v^2+\sigma
_\varepsilon^2)/f_{{\bolds\theta}_0}(u)$. This result is similar to Theorem 1(i) in \citet{KZ12}.

For the dense longitudinal data case, $m_ih\nu_0f_{{\bolds\theta
}_0}(u)(\sigma_v^2+\sigma_\varepsilon^2)$ is dominated by
$m_i(m_i-1)h^2\mu_0^2f_{{\bolds\theta}_0}^2(u)\sigma_v^2$ if
we assume that $m_ih\rightarrow\infty$. Then, again by Lemma 1 in the
supplementary material [\citet{CLLWsupp}], we can prove that
\begin{eqnarray}
\sum_{i=1}^n\mathrm{E} \bigl[{\mathbf
s}_i(u\mid{\bolds\theta}_0){\mathbf
e}_i{\mathbf e}_i^\top{\mathbf
s}_i^\top(u\mid{\bolds\theta}_0)
\bigr]&\sim&\frac{1}{(nh)^2}\sum_{i=1}^n
\frac
{m_i(m_i-1)h^2\mu_0^2\sigma
_v^2}{m_i^2}
\nonumber
\\
&\sim&\frac{\mu_0^2\sigma
_v^2}{n}.
\nonumber
\end{eqnarray}
Hence, in this case, we can take $\varphi_n(h)=n$ and $\sigma_*^2=\mu
_0^2\sigma_v^2$, which are analogous
to those in Theorem 1(ii) of \citet{KZ12} and quite different
from those in the sparse longitudinal data case.
}
\end{remark}


\section{Estimation of covariance matrices}
\label{sec4}

Estimation of the weight or working covariance matrices, which are
involved in the SGEE (\ref{eq27}), is critical to improving the
efficiency of the proposed semiparametric estimators. However, the
unbalanced longitudinal data structure, which can be either sparse or
dense, makes such covariance matrix estimation very challenging, and
some existing estimation methods based on balanced data [such as \citet{W11}] cannot be directly used here. In this section, we introduce a
semiparametric estimation approach that is applicable to both sparse
and dense unbalanced longitudinal data. This approach is based on a
variance--correlation decomposition, and the estimation of the working
covariance matrices then consists of two steps: first, estimate the
conditional variance function using a robust nonparametric method that
accommodates heavy-tailed errors, and second, estimate the parameters
in the correlation matrix. For recent developments on the study of the
covariance structure in longitudinal data analysis, we refer to \citet{FW08},
\citet{ZLT} and the references therein.

For each $1\leq i\leq n$, let ${\mathbf{R}}_i$ be the covariance
matrix of ${\mathbf e}_i$ and
\[
{\bolds\Sigma}_i=\operatorname{diag} \bigl\{\sigma^2(t_{i1}),
\dots,\sigma^2(t_{im_i}) \bigr\}
\]
with $\sigma^2(t_{ij})=\mathrm{E} [e_i^2(t_{ij})\mid t_{ij}
]=\mathrm{E} [e_i^2(t_{ij})\mid t_{ij},{\mathbf X}_i(t_{ij}),
{\mathbf Z}_i(t_{ij}) ]$ for $j=1,\dots,m_i$, and ${\mathbf
{C}}_i$ be the correlation matrix of ${\mathbf e}_i$. Assume that there
exists a $q$-dimensional parameter vector ${\bolds\phi}$ such
that ${\mathbf C}_i={\mathbf{C}}_i({\bolds\phi})$ where
${\mathbf{C}}_i(\cdot)$, $1\leq i\leq n$, are pre-specified. By the
variance--correlation decomposition, we have
%
\begin{equation}
\label{eq41} {\mathbf{R}}_i={\bolds\Sigma}_i^{1/2}{
\mathbf{C}}_i({\bolds\phi}){\bolds\Sigma}_i^{1/2}.
\end{equation}
The above semiparametric covariance structure has been studied in some
of the existing literature [see, e.g., \citet{FHL07} and \citet{FW08}] and provides a flexible framework to capture the error
covariance structure, especially when the dimension of ${\bolds\phi}$ is large. For example, it is satisfied when $e_i(t_{ij})$ has
the $\operatorname{AR}(1)$ or $\operatorname{ARMA}(1,1)$ dependence structure for each $i$; see, for
example, the simulated example in Section~\ref{sec5.1}. When\vspace*{2pt}
$e_{i}(t_{ij})=\sigma(t_{ij})(v_i+\varepsilon_{ij})$ in which $v_i$
and $\varepsilon_{ij}$ satisfy the conditions discussed in Remark~\ref{re3}
and $\sigma^2_{\varepsilon}+\sigma^2_v=1$, we can also show that the
semiparametric covariance structure is satisfied with ${\bolds\phi}$ being $\sigma^2_\varepsilon$ or $\sigma^2_v$. Some existing
papers such as \citet{WP03} suggest the use of a
nonparametric smoothing method to estimate the covariance matrix.
However, they usually need to assume that the longitudinal data are
balanced or nearly balanced, which would be violated when the data are
collected at irregular and possibly subject-specific time points.
\citet{YMW} proposed the approach of functional data analysis to
estimate the covariance structure for sparse and irregularly-spaced
longitudinal data. However, some substantial modification may be needed
to extend the method of \citet{YMW} to our framework, which
includes both the sparse and dense longitudinal data.

In the present paper, we first estimate the conditional variance function
$\sigma^2(\cdot)$ in the diagonal matrix ${\bolds\Sigma}_i$ by
using a nonparametric method. In recent years, there has been a rich
literature on the study of nonparametric conditional variance
estimation; see, for example, \citet{FY98}, \citet{YM04},
\citet{FHL07} and \citet{LT11}. However, when the errors are heavy-tailed, which is not
uncommon in
economic and financial data analysis, most of these existing methods
may not perform well. This motivates us to devise an estimation
method that is robust to heavy-tailed errors. Let
$r(t_{ij})= [Y_{ij}-{\mathbf{Z}}_{ij}^\top{\bolds\beta}_0
-\eta({\mathbf{X}}_{ij}^\top{\bolds\theta}_0) ]^2$. We
can then find a
random variable $\xi(t_{ij})$ so that
$r(t_{ij})=\sigma^2(t_{ij})\xi^2(t_{ij})$ and $\mathrm{E} [\xi^2(t_{ij})\mid t_{ij} ]=1$ with probability one. By
applying the log-transformation [see \citet{PY03} and \citet{CCP09} for the application of this
transformation in time series analysis] to $r(t_{ij})$, we have
%
\begin{equation}
\label{eq42} \log r(t_{ij})=\log\bigl[\tau\sigma^2(t_{ij})
\bigr]+\log\bigl[\tau^{-1}\xi^2(t_{ij}) \bigr]
\equiv\sigma^2_\diamond(t_{ij})+\xi
_{\diamond}(t_{ij}),
\end{equation}
where $\tau$ is a positive constant such that $\mathrm{E}[\xi_{\diamond
}(t_{ij})]=\mathrm{E} \{\log[\tau^{-1}\xi^2(t_{ij}) ]
\}=0$. Here, $\xi_{\diamond}(t_{ij})$ could be viewed as an error
term in model (\ref{eq42}). As $r_{ij}\equiv r(t_{ij})$ are
unobservable, we replace them with
\[
\widehat{r}_{ij}= \bigl[Y_{ij}-{\mathbf{Z}}_{ij}^\top
\widetilde{\bolds\beta}-\widehat{\eta}\bigl({\mathbf{X}}_{ij}^\top
\widetilde{\bolds\theta}\mid\widetilde{\bolds\beta},
\widetilde{
\bolds\theta}\bigr) \bigr]^2,
\]
where $\widetilde{\bolds\beta}$ and $\widetilde{\bolds\theta}$ are the PULS estimators of ${\bolds\beta}_0$ and
${\bolds\theta}_0$, respectively. In order to estimate $\sigma
^2_\diamond(t)$, we define
%
\begin{equation}
\label{eq43} \qquad\widetilde{L}_n (a,b )=\sum
_{i=1}^n \Biggl\{\frac
{w_i}{h_1}\sum
_{j=1}^{m_i} \bigl[\log(\widehat{r}_{ij}+\zeta
_n)-a-b(t_{ij}-t) \bigr]^2 K_1
\biggl(\frac{t_{ij}-t}{h_1} \biggr) \Biggr\},
\end{equation}
where $K_1(\cdot)$ is a kernel function, $h_1$ is a bandwidth
satisfying Assumption~\ref{as9} in Appendix~\ref{appA}, $w_i={1}/{(nm_i)}$ as in
Section~\ref{sec2} and $\zeta_n\rightarrow0$ as $n\rightarrow\infty$.
Throughout this paper, we set $\zeta_n={1}/{T_n}$, where
$T_n=\sum_{i=1}^n m_i$. The $\zeta_n$ is added in $\log
(\widehat{r}_{ij}+\zeta_n)$ to avoid the occurrence of invalid $\log
0$ as $\zeta_n>0$ for any $n$. Such a modification would not affect
the asymptotic distribution of the conditional variance estimation
under certain mild restrictions. Then $\sigma^2_\diamond(t)$ can be
estimated as
%
\begin{equation}
\label{eq44} \widehat{\sigma}_\diamond^2(t)=\widehat{a}\qquad\mbox{where } (\widehat{a},\widehat{b})^\top=\arg\min
_{a,b}\widetilde{L}_n(a,b).
\end{equation}
On the other hand, noting that $\frac{\exp\{\sigma^2_\diamond
(t_{ij}) \}}{\tau}\xi^2(t_{ij})=r_{ij}$ and $\mathrm{E}[\xi
^2(t_{ij})]=1$, the constant $\tau$ can be estimated by
%
\begin{equation}
\label{eq45} \widehat{\tau}= \Biggl[\frac{1}{T_n}\sum
_{i=1}^n\sum_{j=1}^{m_i}
\widehat{r}_{ij}\exp\bigl\{-\widehat{\sigma}^2_\diamond
(t_{ij})\bigr\} \Biggr]^{-1}.
\end{equation}
We then estimate $\sigma^2(t)$ by
%
\begin{equation}
\label{eq46} \widehat{\sigma}^2(t)=\frac{\exp\{\widehat{\sigma
}^2_\diamond
(t)\}}{\widehat{\tau}}.
\end{equation}
It is easy to see that thus defined estimator $\widehat{\sigma}^2(t)$
is always positive.

Suppose that there exists a sequence $\varphi_{n\diamond}(h_1)$ which
depends on $h_1$, and a constant $0<\sigma^2_\diamond<\infty$ such that
%
\begin{eqnarray}\label{eq47}
\varphi_{n\diamond}(h_1)&=&o(\omega_n),
\nonumber\\[-8pt]\\[-8pt]\nonumber
\frac{\varphi_{n\diamond
}(h_1)}{h_1^2}\max_{1\leq i\leq n}w_i^2\mathrm{E} \Biggl[\sum_{j=1}^{m_i}
\xi_\diamond(t_{ij})K_1 \biggl(\frac{t_{ij}-t}{h_1}
\biggr) \Biggr]^2&=&o(1)
\end{eqnarray}
and
%
\begin{equation}
\label{eq471} \frac{\varphi_{n\diamond}(h_1)}{f_T(t)h_1^2}\mathrm{E} \Biggl
[\sum
_{i=1}^nw_i\sum
_{j=1}^{m_i}\xi_\diamond(t_{ij})K_1
\biggl(\frac
{t_{ij}-t}{h_1} \biggr) \Biggr]^2\rightarrow
\sigma^2_\diamond,
\end{equation}
which are similar to those in (\ref{eq340}) and (\ref{eq34}),
where $f_T(\cdot)$ is the density function of the observation times
$t_{ij}$. Define
\begin{eqnarray}
b_{\sigma1}(t)&=&\frac{\exp\{\sigma^2_\diamond(t)\}}{2\tau}\ddot{\sigma
}_\diamond^2(t)
\int v^2K_1(v)\,dv,
\nonumber
\\
b_{\sigma2}(t)&=&\frac{\exp\{\sigma^2_\diamond(t)\}}{2\tau}\mathrm{E}\bigl
[\ddot{\sigma}_\diamond^2(t_{ij})
\bigr]\int v^2K_1(v)\,dv,
\nonumber
\end{eqnarray}
where $\ddot{\sigma}^2_\diamond(\cdot)$ is the second-order
derivative of $\sigma^2_\diamond(\cdot)$. We then establish the
asymptotic distribution of $\widehat{\sigma}^2(t)$ in the following
theorem, whose proof is given in the supplementary material [\citet{CLLWsupp}].

\begin{theorem}\label{teo3}
Suppose the conditions in Theorems~\ref{teo1} and~\ref{teo2}, Assumptions~\ref{as6}--\ref{as9} in
Appendix~\ref{appA}, (\ref{eq47}) and (\ref{eq471}) are satisfied. Then we have
%
\begin{equation}
\label{eq48} \qquad\varphi_{n\diamond}^{1/2}(h_1) \bigl\{
\widehat{\sigma}^2(t)-\sigma^2(t)-\bigl[b_{\sigma1}(t)-b_{\sigma2}(t)
\bigr]h_1^2 \bigr\}\stackrel{d}\longrightarrow\mathrm{N}
\biggl(0,\frac{\sigma^4(t)}{f_T(t)}\sigma^2_\diamond\biggr).
\end{equation}
\end{theorem}

\begin{remark}\label{re4}
{Theorem~\ref{teo3} can be seen as an extension of Theorem 1 in \citet{CCP09}
from the time series case to the longitudinal data case. The
longitudinal data framework in this paper is more flexible and includes
both sparse and dense data types. If $\xi_{\diamond
}(t_{ij})=v_i^{\diamond}+\varepsilon_{ij}^\diamond$, where
$\varepsilon_{ij}^\diamond$ are i.i.d. across both $i$ and $j$ with
$\mathrm{E}[\varepsilon_{ij}^\diamond]=0$ and $\mathrm{E}[(\varepsilon
_{ij}^\diamond)^2]<\infty$, and $\{v_{i}^\diamond\}$ is an i.i.d.
sequence of random variables with $\mathrm{E}[v_{i}^\diamond]=0$ and
$\mathrm{E}[(v_{i}^\diamond)^2]<\infty$ and is independent of $\{
\varepsilon_{ij}^\diamond\}$,\vspace*{2pt} following the discussion in Remark~\ref{re3},
we can again show that the form of $\varphi_{n\diamond}(h_1)$ depends
on the type of the longitudinal data, and thus the nonparametric
conditional variance estimation has different convergence rates for
sparse and dense data.
}
\end{remark}

We next\vspace*{1pt} discuss how to obtain the optimal value of the parameter vector
${\bolds\phi}$. Construct the residuals $\widetilde{\mathbf
e}_i={\mathbf Y}_i-{\mathbf Z}_i\widetilde{\bolds\beta
}-\widetilde{\bolds\eta}({\mathbf X}_i,\widetilde{\bolds\theta})$, where $\widetilde{\bolds\eta}({\mathbf
X}_i,\widetilde{\bolds\theta})$ is defined in the same\vspace*{1pt} way as
${\bolds\eta}({\mathbf X}_i,{\bolds\theta})$ but with
$\eta(\cdot)$ and ${\bolds\theta}$ replaced by $\widetilde
{\eta}(\cdot)\equiv\widehat{\eta}(\cdot\mid\widetilde{\bolds\beta},\widetilde{\bolds\theta})$ and $\widetilde{\bolds\theta}$, respectively. Let $\widetilde{\bolds\Lambda}_i\equiv
\widehat{\bolds\Lambda}_i(\widetilde{\bolds\theta})$,
$\widehat{\bolds\Sigma}_i=\operatorname{diag} \{\widehat{\sigma
}^2(t_{i1}),\dots,\widehat{\sigma}^2(t_{im_i}) \}$, and define
${\mathbf{R}}_i^*({\bolds\phi})=\widehat{\bolds\Sigma
}_i^{1/2}{\mathbf{C}}_i({\bolds\phi})\widehat{\bolds\Sigma}_i^{1/2}$. Motivated by equations (\ref{eq31}) and (\ref
{eq32}), we construct
%
\begin{equation}
\label{eq49} {\bolds\Omega}_0^*({\bolds\phi})=\sum
_{i=1}^n\widetilde{\bolds\Lambda}_i^\top\bigl[{\mathbf{R}}_i^{*}({
\bolds\phi}) \bigr]^{-1}\widetilde{\bolds\Lambda}_i
\end{equation}
and
%
\begin{equation}
\label{eq49a} {\bolds\Omega}_1^*({\bolds\phi})=\sum
_{i=1}^n\widetilde{\bolds\Lambda}_i^\top\bigl[{\mathbf{R}}_i^{*}({
\bolds\phi}) \bigr]^{-1}\widetilde{\mathbf e}_i
\widetilde{\mathbf e}_i^\top\bigl[{\mathbf{R}}_i^{*}({
\bolds\phi}) \bigr]^{-1}\widetilde{\bolds\Lambda}_i.
\end{equation}
By Theorem~\ref{teo1}, the sandwich formula estimate $ [{\bolds\Omega
}_0^{*}({\bolds\phi}) ]^+{\bolds\Omega
}_1^*({\bolds\phi}) [{\bolds\Omega}_0^*({\bolds\phi}) ]^+$ is asymptotically proportional to the asymptotic
covariance of the proposed SGEE estimators when the inverse of
${\mathbf{R}}_i^*({\bolds\phi})$ is chosen as the weight
matrix. The optimal value of ${\bolds\phi}$, denoted by
$\widehat{\bolds\phi}$, can be chosen to minimize the
determinant $\mid [{\bolds\Omega}_0^{*}({\bolds\phi
}) ]^+{\bolds\Omega}_1^*({\bolds\phi})
[{\bolds\Omega}_0^*({\bolds\phi}) ]^+\mid$. Such a
method is called the minimum generalized variance method [\citet{FHL07}]. With the chosen $\widehat{\bolds\phi}$, we can
estimate the covariance matrices by
%
\begin{equation}
\label{eq491} {\mathbf{R}}_i(\widehat{\bolds\phi})=\widehat{
\bolds\Sigma}_i^{1/2}{\mathbf{C}}_i(
\widehat{\bolds\phi})\widehat{\bolds\Sigma}_i^{1/2},
\end{equation}
whose inverse will be used as the weight matrices in the SGEE method.


\section{Numerical studies}\label{sec5}

In this section, we first study the finite sample performance of the
proposed SGEE estimators through Monte Carlo simulation, and then give
an empirical application of the proposed model and methodology.

\subsection{Simulation studies}\label{subsec51}\label{sec5.1}


We investigate both sparse and dense longitudinal data cases with an
average time dimension $\overline{m}$ of $10$ for the sparse data and
$30$ for the dense data. We use two types of within-subject correlation
structure, $\operatorname{AR}(1)$ and $\operatorname{ARMA}(1,1)$, in the error terms $e_i(t_{ij})$. We
investigate the finite sample performance of the proposed estimators
under both correct specification and misspecification of the
correlation structure in the construction of the covariance matrix
estimator proposed in Section~\ref{sec4}. For the misspecified case, we fit an
$\operatorname{AR}(1)$ correlation structure while the true underlying structure is
$\operatorname{ARMA}(1,1)$ and examine the robustness of the estimators.


Simulated data are generated from model (\ref{eq12}) with
two-dimensional ${\mathbf Z}_i(t_{ij})$ and three-dimensional ${\mathbf
X}_i(t_{ij})$, and
\[
{\bolds\beta}_0=(2,1)^\top,\qquad {\bolds\theta
}_0=(2,1,2)^\top/3\quad\mbox{and}\quad \eta(u)=0.5\exp(u).
\]
The covariates $({\mathbf Z}_i^\top(t_{ij}), {\mathbf X}_i^\top
(t_{ij}))^\top$ are generated independently from a five-dimensional
Gaussian distribution with mean ${\mathbf0}$, variance 1 and pairwise
correlation $0.1$. 
The observation times $t_{ij}$ are generated in the same way as in
\citet{FHL07}: for each subject, $\{0,1,2,\dots, T\}$ is a set of
scheduled times, and each scheduled time from 1 to $T$ has a $0.2$
probability of being skipped; each actual observation time is a
perturbation of a nonskipped scheduled time; that is, a uniform $[0,
1]$ random number is added to the nonskipped scheduled time. Here $T$
is set to be 12 or 36, which corresponds to an average time dimension
of $\overline{m}=10$ or $\overline{m}=30$, respectively. For each
$i$, the error terms $e_i(t_{ij})$ are generated from a Gaussian
process with mean 0, variance function
%
\begin{equation}
\label{eq51} \operatorname{var}\bigl[e(t)\bigr]=\sigma^2(t)=0.25\exp(t/12)
\end{equation}
and serial correlation structure
%
\begin{equation}
\label{eq52} \operatorname{cor}\bigl(e(t),e(s)\bigr)=\cases{ 1, &\quad$t=s$,
\vspace*{3pt}\cr
\gamma
\rho^{|t-s|}, &\quad$t\neq s$.}
\end{equation}
Note that (\ref{eq52}) corresponds to an $\operatorname{ARMA}(1,1)$ correlation
structure and reduces to an $\operatorname{AR}(1)$ correlation structure when $\gamma
=1$. The number of subjects, $n$, is taken to be $30$ or $50$. The
values for $\gamma$ and $\rho$ are $(\gamma,\rho)=(0.85,0.9)$ in
the $\operatorname{ARMA}(1,1)$ correlation structure and $(\gamma,\rho)=(1,0.9)$ in
the $\operatorname{AR}(1)$ structure.

\begin{table}
\tabcolsep=0pt
\caption{Performance of parameter estimation methods under correct
specification of an underlying $\operatorname{AR}(1)$ correlation structure}\label{table1}
\begin{tabular*}{\tablewidth}{@{\extracolsep{\fill}}@{}lccd{2.4}ccd{2.4}cc@{}}
\hline
& \multicolumn{2}{c}{$\bolds{n}$} & \multicolumn{3}{c}{\textbf{30}} & \multicolumn{3}{c@{}}{\textbf{50}}\\[-6pt]
& \multicolumn{2}{c}{\hrulefill} & \multicolumn{3}{c}{\hrulefill} & \multicolumn{3}{c@{}}{\hrulefill}\\
$\bolds{\overline{m}}$ & \multicolumn{1}{c}{\textbf{Parameters}} & \multicolumn{1}{c}{\textbf{Methods}} &
\multicolumn{1}{c}{\textbf{Bias}} & \multicolumn{1}{c}{\textbf{SD}} & \multicolumn{1}{c}{\textbf{MAD}} &
\multicolumn{1}{c}{\textbf{Bias}} & \multicolumn{1}{c}{\textbf{SD}} & \multicolumn{1}{c@{}}{\textbf{MAD}}\\
\hline
10 & {$\beta_1$} & PULS & 0.0048 & 0.0402 & 0.0288 & -0.0030 & 0.0308 & 0.0195 \\
 & & SGEE & -0.0026 & 0.0508 &0.0081 & -0.0016 & 0.0259 & 0.0074\\[3pt]
 &  {$\beta_2$} & PULS & -0.0024 & 0.0409 & 0.0243 & 0.0049 & 0.0267 & 0.0180 \\
 & & SGEE & -0.0018 & 0.0298 & 0.0110 & 0.0033 & 0.0310 & 0.0077\\[3pt]
  &  {$\theta_1$} & PULS & -0.0049 & 0.0299 & 0.0180 & -0.0009 & 0.0197 & 0.0134 \\
 & & SGEE & -0.0013 & 0.0164 & 0.0083 & -0.0002 & 0.0118 & 0.0046\\[3pt]
  &  {$\theta_2$} & PULS & 0.0011 & 0.0380 & 0.0229 & -0.0016 & 0.0237 & 0.0161 \\
 & & SGEE & 0.0026 & 0.0188 &0.0100 & 0.0006 & 0.0108 & 0.0067\\[3pt]
  &  {$\theta_3$} & PULS & 0.0018 & 0.0314 & 0.0188 & 0.0006 & 0.0203 & 0.0147 \\
 & & SGEE & -0.0007 & 0.0182 & 0.0090 & -0.0004 & 0.0088 & 0.0052
 \\[6pt]
 {30} &  {$\beta_1$} & PULS & 0.0003 & 0.0408 & 0.0277 & 0.0016 & 0.0328 & 0.0222 \\
 & & SGEE & -0.0081 & 0.1134 &0.0106 & 0.0007 & 0.0108 & 0.0083\\[3pt]
  &  {$\beta_2$} & PULS & -0.0020 &0.0425 & 0.0317 & 0.0005 & 0.0351 & 0.0202 \\
 & & SGEE & -0.0017 & 0.0420 & 0.0096 & -0.0064 &0.0152 &0.0079\\[3pt]
  &  {$\theta_1$} & PULS & 0.0020 &0.0315 & 0.0213 & -0.0020 & 0.0244 & 0.0182 \\[3pt]
 & & SGEE & -0.0008 & 0.0247 & 0.0075 & 0.0001 &0.0148 &0.0064\\[3pt]
  &  {$\theta_2$} & PULS & -0.0035 &0.0340 & 0.0240 & -0.0083 & 0.0278 & 0.0163 \\
 & & SGEE & -0.0027 & 0.0242 &0.0090 & -0.0013 &0.0104 &0.0066\\[3pt]
  &  {$\theta_3$} & PULS & -0.0027 &0.0321 & 0.0185 & 0.0045 & 0.0267 & 0.0169 \\
 & & SGEE & 0.0009 & 0.0230 & 0.0074 & 0.0001 &0.0162 &0.0068\\
\hline
\end{tabular*}
\end{table}

\begin{figure}
\begin{tabular}{cc}

\includegraphics{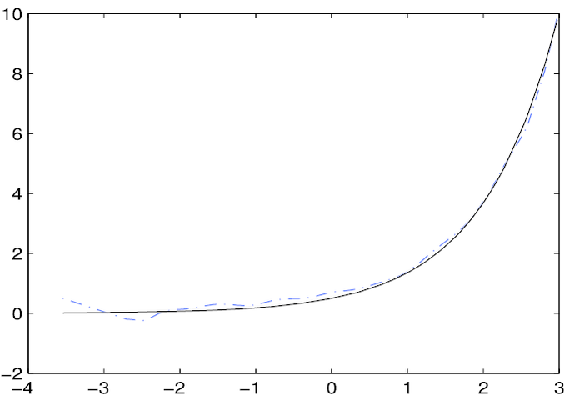}
 & \includegraphics{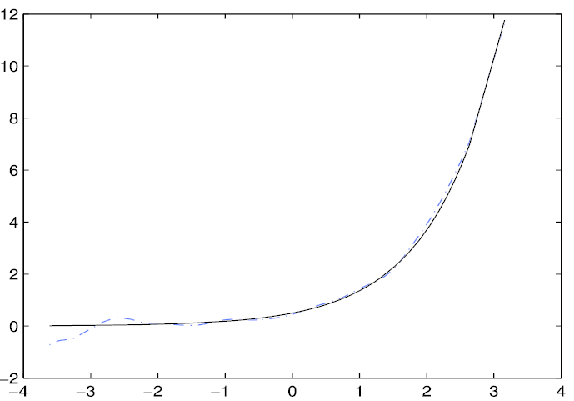}\\
{\footnotesize(a)} & {\footnotesize(b)}
\\[6pt]

\includegraphics{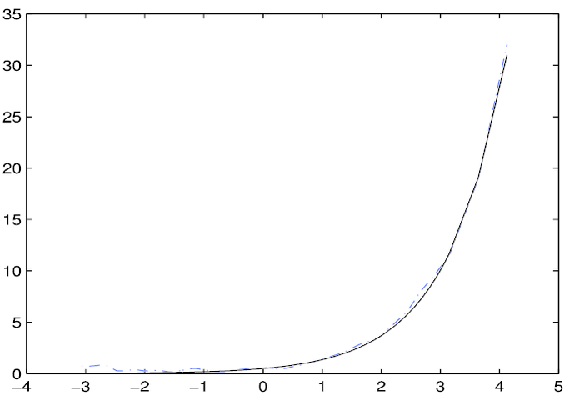}
 & \includegraphics{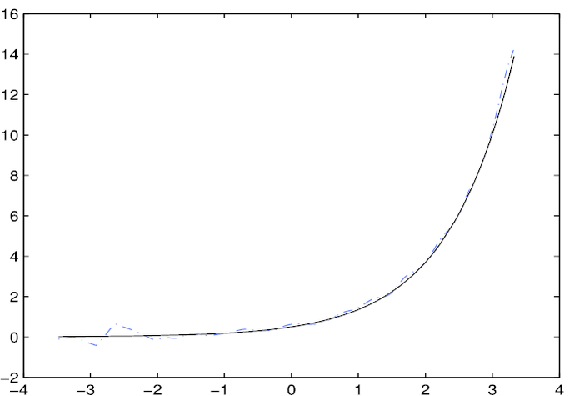}\\
{\footnotesize(c)} & {\footnotesize(d)}
\end{tabular}
\caption{Estimated link function (dot-dashed line), together with the
true link function (solid line), from a typical realization of model
(\protect\ref{eq12}) with $\operatorname{AR}(1)$ correlation structure for each
combination of $n$ and $\overline{m}$:
\textup{(a)} $n=30$, $\overline{m}=10$;
\textup{(b)} $n=50$, $\overline{m}=10$;
\textup{(c)} $n=30$, $\overline{m}=30$;
\textup{(d)} $n=50$, $\overline{m}=30$.}\label{fig1}
\end{figure}

\begin{figure}
\begin{tabular}{cc}

\includegraphics{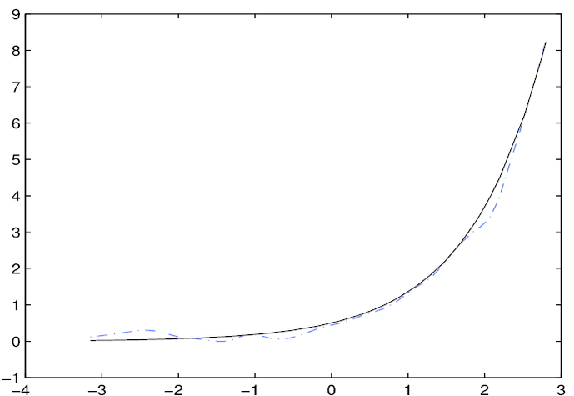}
 & \includegraphics{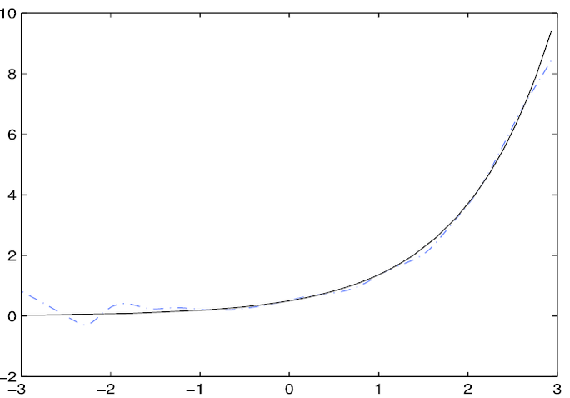}\\
{\footnotesize(a)} & {\footnotesize(b)}
\\[6pt]

\includegraphics{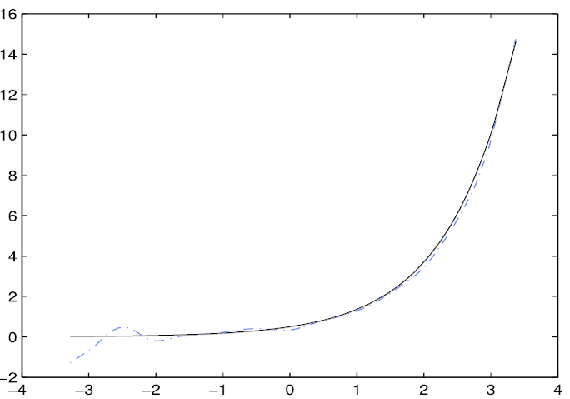}
 & \includegraphics{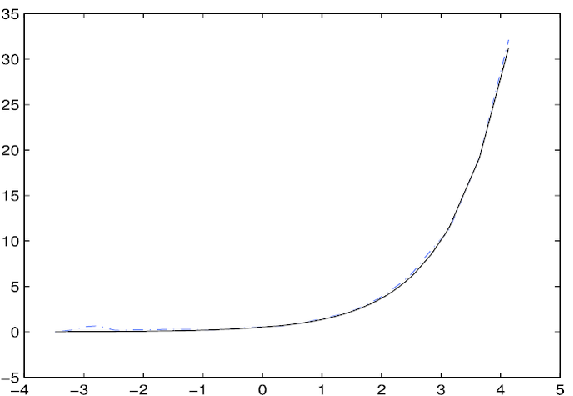}\\
{\footnotesize(c)} & {\footnotesize(d)}
\end{tabular}
\caption{Estimated link function (dot-dashed line), together with the
true link function (solid line), from a typical realization of model
(\protect\ref{eq12}) with $\operatorname{ARMA}(1,1)$ correlation structure for each
combination of $n$ and $\overline{m}$:
\textup{(a)} $n=30$, $\overline{m}=10$;
\textup{(b)} $n=50$, $\overline{m}=10$;
\textup{(c)} $n=30$, $\overline{m}=30$;
\textup{(d)} $n=50$, $\overline{m}=30$.}\label{fig2}
\end{figure}

For each combination of $\overline{m}$, $n$, and the correlation
structure, the number of simulation replications is 200. For the
selection of the bandwidth, however, due to the running time
limitation, we first run a leave-one-unit-out (i.e., leave out
observations from one subject at a time) cross-validation (CV) to
choose the optimal bandwidths from 20 replications. We then use the
average of the optimal bandwidths from these 20 replications as the
bandwidth in the 200 replications of the simulation study. For the SGEE
method, we choose the weight matrix as the inverse of the estimated
within-subject covariance matrix as constructed in (\ref{eq491}) of
Section~\ref{sec4}. We first study the performance of the proposed estimators in
the case where the correlation structure in the estimation of the
covariance matrix is correctly specified, and then investigate the
robustness of the estimators to the misspecification of the correlation
structure. The bias, calculated as the average of the estimates from
the 200 replications minus the true parameter values, the standard
deviation (SD), calculated as the sample standard deviation of the 200
estimates and the median absolute deviation (MAD), calculated as the
median absolute deviation of the 200 estimates are reported in Tables~\ref{table1}
and \ref{table2}. Table~\ref{table1} gives the results obtained under the correct
specification of an underlying within-subject $\operatorname{AR}(1)$ correlation
structure in $e_i(t_{ij})$, and Table~\ref{table2} gives those obtained under the
correct specification of an underlying $\operatorname{ARMA}(1,1)$ structure in
$e_i(t_{ij})$. For comparison, we also report the results from the PULS
estimation. The results in Tables~\ref{table1} and \ref{table2} show that the SGEE estimates
are comparable with the PULS estimates in terms of bias and are more
efficient than the PULS estimates, which supports the asymptotic theory
developed in Section~\ref{sec3}. In Figures~\ref{fig1} and \ref{fig2}, we plot the local linear
estimated link function from a typical realization together with the
real curve for each combination of $n$ and $\overline{m}$.


\begin{table}
\caption{Performance of parameter estimation methods under correct
specification of an underlying $\operatorname{ARMA}(1,1)$ correlation structure}\label{table2}
\begin{tabular*}{\tablewidth}{@{\extracolsep{\fill}}@{}lccd{2.4}ccd{2.4}cc@{}}
\hline
& \multicolumn{2}{c}{$\bolds{n}$} & \multicolumn{3}{c}{\textbf{30}} & \multicolumn{3}{c@{}}{\textbf{50}}\\[-6pt]
& \multicolumn{2}{c}{\hrulefill} & \multicolumn{3}{c}{\hrulefill} & \multicolumn{3}{c@{}}{\hrulefill}\\
$\bolds{\overline{m}}$ & \multicolumn{1}{c}{\textbf{Parameters}} & \multicolumn{1}{c}{\textbf{Methods}} &
\multicolumn{1}{c}{\textbf{Bias}} & \multicolumn{1}{c}{\textbf{SD}} & \multicolumn{1}{c}{\textbf{MAD}} &
\multicolumn{1}{c}{\textbf{Bias}} & \multicolumn{1}{c}{\textbf{SD}} & \multicolumn{1}{c@{}}{\textbf{MAD}}\\
\hline
{10} &  {$\beta_1$} & PULS & -0.0029 & 0.0400 & 0.0280 & 0.0006 & 0.0322 & 0.0221 \\
 & & SGEE & -0.0025 & 0.0244 &0.0155 & 0.0000 & 0.0193 & 0.0124\\[3pt]
  &  {$\beta_2$} & PULS & 0.0032 & 0.0386 & 0.0282 & -0.0045 & 0.0299 & 0.0205 \\
 & & SGEE & 0.0009 & 0.0249 & 0.0171 & 0.0001 & 0.0212 & 0.0126\\[3pt]
  &  {$\theta_1$} & PULS & -0.0004 & 0.0267 & 0.0181 & -0.0003 & 0.0188 & 0.0126 \\
 & & SGEE & -0.0002 & 0.0161 & 0.0104 & 0.0006 & 0.0146 & 0.0073\\[3pt]
  &  {$\theta_2$} & PULS & -0.0047 & 0.0343 & 0.0209 & 0.0005 & 0.0223 & 0.0156 \\
 & & SGEE & -0.0031 & 0.0192 &0.0113 & -0.0002 & 0.0145 & 0.0087\\[3pt]
  &  {$\theta_3$} & PULS & 0.0008 & 0.0253 & 0.0158 & -0.0009 & 0.0201 & 0.0121 \\
 & & SGEE & 0.0011 & 0.0148 & 0.0102 & -0.0009 & 0.0146 & 0.0074\\[6pt]
{30} &  {$\beta_1$} & PULS & -0.0026 & 0.0450 & 0.0296 & -0.0016 & 0.0374 & 0.0273 \\
 & & SGEE & 0.0005 & 0.0214 &0.0138 & 0.0015 & 0.0288 & 0.0105\\[3pt]
  &  {$\beta_2$} & PULS & -0.0013 & 0.0461 & 0.0291 & 0.0035 & 0.0361 & 0.0252 \\
 & & SGEE & 0.0040 & 0.0335 & 0.0147 & 0.0014 & 0.0152 &0.0104\\[3pt]
  &  {$\theta_1$} & PULS & -0.0014 & 0.0296 & 0.0192 & -0.0010 & 0.0207 & 0.0159 \\
 & & SGEE & -0.0005 & 0.0166 & 0.0095 & 0.0006 &0.0092 &0.0063\\[3pt]
  &  {$\theta_2$} & PULS & -0.0050 &0.0355 & 0.0231 & 0.0011 & 0.0229 & 0.0173 \\
 & & SGEE & -0.0037 & 0.0371 &0.0120 & -0.0003 &0.0116 &0.0072\\[3pt]
  &  {$\theta_3$} & PULS & 0.0017 & 0.0279 & 0.0186 & -0.0006 & 0.0215 & 0.0154 \\
 & & SGEE & 0.0009 & 0.0181 & 0.0095 & -0.0007 &0.0100 &0.0070\\
\hline
\end{tabular*}
\end{table}


\begin{table}
\caption{Performance of parameter estimation methods under
misspecification of an underlying $\operatorname{ARMA}(1,1)$ correlation structure}\label{table3}
\begin{tabular*}{\tablewidth}{@{\extracolsep{\fill}}@{}lccd{2.4}ccd{2.4}cc@{}}
\hline
& \multicolumn{2}{c}{$\bolds{n}$} & \multicolumn{3}{c}{\textbf{30}} & \multicolumn{3}{c@{}}{\textbf{50}}\\[-6pt]
& \multicolumn{2}{c}{\hrulefill} & \multicolumn{3}{c}{\hrulefill} & \multicolumn{3}{c@{}}{\hrulefill}\\
$\bolds{\overline{m}}$ & \multicolumn{1}{c}{\textbf{Parameters}} & \multicolumn{1}{c}{\textbf{Methods}} &
\multicolumn{1}{c}{\textbf{Bias}} & \multicolumn{1}{c}{\textbf{SD}} & \multicolumn{1}{c}{\textbf{MAD}} &
\multicolumn{1}{c}{\textbf{Bias}} & \multicolumn{1}{c}{\textbf{SD}} & \multicolumn{1}{c@{}}{\textbf{MAD}}\\
\hline
{10} &  {$\beta_1$} & PULS &0.0072 & 0.0410 & 0.0357 & -0.0038 & 0.0299 & 0.0201 \\
 & & SGEE & -0.0054 & 0.0261 &0.0210 & -0.0055 &0.0211 & 0.0147\\[3pt]
  &  {$\beta_2$} & PULS & 0.0068 &0.0336 & 0.0256 & 0.0037 & 0.0290 & 0.0163 \\
 & & SGEE & 0.0025 & 0.0267 & 0.0157 & 0.0023 &0.0190 & 0.0136\\[3pt]
  &  {$\theta_1$} & PULS & 0.0037 &0.0166 & 0.0114 & 0.0061 & 0.0157 & 0.0096 \\
 & & SGEE & 0.0033 & 0.0144 & 0.0122 & 0.0016 &0.0163 & 0.0081\\[3pt]
  &  {$\theta_2$} & PULS & -0.0092 &0.0303 & 0.0184 & -0.0084 & 0.0224 & 0.0174 \\
 & & SGEE & -0.0007 & 0.0198 &0.0144 & -0.0045 &0.0203 & 0.0130\\[3pt]
  &  {$\theta_3$} & PULS & -0.0005 &0.0229 & 0.0158 & -0.0028 & 0.0160 & 0.0111 \\
 & & SGEE & -0.0035 & 0.0141 & 0.0094 & 0.0000 &0.0134 & 0.0092\\[6pt]
{30} &  {$\beta_1$} & PULS &0.0066 & 0.0403 & 0.0259 & -0.0221 & 0.0502 & 0.0252 \\
 & & SGEE & 0.0093 & 0.0144 &0.0087 &0.0001 & 0.0165& 0.0118\\[3pt]
  &  {$\beta_2$} & PULS & -0.0138 &0.0435 & 0.0353 & 0.0107 & 0.0312 & 0.0233 \\
 & & SGEE & -0.0017 & 0.0268 & 0.0096 & 0.0035 &0.0170 &0.0096\\[3pt]
  &  {$\theta_1$} & PULS & 0.0027 &0.0252 & 0.0165 & 0.0020 & 0.0181 & 0.0067 \\
 & & SGEE & 0.0054 & 0.0136 & 0.0078 & 0.0019 &0.0096&0.0098\\[3pt]
  &  {$\theta_2$} & PULS & -0.0063 &0.0265 & 0.0245 & 0.0021 & 0.0315 & 0.0273 \\
 & & SGEE & 0.0009 & 0.0198 &0.0118 & 0.0046 & 0.0136&0.0094\\[3pt]
  &  {$\theta_3$} & PULS & -0.0011 &0.0285 & 0.0258 & -0.0042 & 0.0217 & 0.0136 \\
 & & SGEE & -0.0065 & 0.0178 & 0.0137 & -0.0046 &0.0120 &0.0084\\
\hline
\end{tabular*}
\end{table}

\begin{table}[t]
\caption{Performance of parameter estimation methods under correct
specification of an underlying $\operatorname{AR}(1)$ correlation structure when the
covariates in $\mathbf{Z}$ are discrete}\label{table4}
\begin{tabular*}{\tablewidth}{@{\extracolsep{\fill}}@{}lccd{2.4}ccd{2.4}cc@{}}
\hline
& \multicolumn{2}{c}{$\bolds{n}$} & \multicolumn{3}{c}{\textbf{30}} & \multicolumn{3}{c@{}}{\textbf{50}}\\[-6pt]
& \multicolumn{2}{c}{\hrulefill} & \multicolumn{3}{c}{\hrulefill} & \multicolumn{3}{c@{}}{\hrulefill}\\
$\bolds{\overline{m}}$ & \multicolumn{1}{c}{\textbf{Parameters}} & \multicolumn{1}{c}{\textbf{Methods}} &
\multicolumn{1}{c}{\textbf{Bias}} & \multicolumn{1}{c}{\textbf{SD}} & \multicolumn{1}{c}{\textbf{MAD}} &
\multicolumn{1}{c}{\textbf{Bias}} & \multicolumn{1}{c}{\textbf{SD}} & \multicolumn{1}{c@{}}{\textbf{MAD}}\\
\hline
{10} &  {$\beta_1$} & PULS &0.0215 & 0.0530 & 0.0404 & 0.0018 & 0.0646 & 0.0472 \\
 & & SGEE & 0.0228 & 0.0511 &0.0208 & 0.0037 & 0.0298& 0.0138\\[3pt]
  &  {$\beta_2$} & PULS & -0.0309 &0.0858 & 0.0735 & 0.0193 & 0.0526 & 0.0498 \\
 & & SGEE & 0.0024 & 0.0313 & 0.0193 & 0.0074 &0.0339 & 0.0274\\[3pt]
  &  {$\theta_1$} & PULS & -0.0012 &0.0185 & 0.0090 & -0.0116 & 0.0201 & 0.0175 \\
 & & SGEE & -0.0060 & 0.0157 & 0.0082 & 0.0020 &0.0086 & 0.0066\\[3pt]
  &  {$\theta_2$} & PULS & -0.0020 &0.0263 & 0.0232 & 0.0138 & 0.0229 & 0.0172 \\
 & & SGEE & 0.0122 & 0.0241 & 0.0143 & -0.0004 &0.0087 & 0.0063\\[3pt]
  &  {$\theta_3$} & PULS & 0.0012 &0.0206 & 0.0075 & 0.0036 & 0.0153 & 0.0132 \\
 & & SGEE & -0.0008 & 0.0078 & 0.0048 & -0.0020 &0.0070 & 0.0034\\[6pt]
{30} &  {$\beta_1$} & PULS &0.0075 & 0.0427 & 0.0222 & 0.0108 & 0.0723 & 0.0513 \\
 & & SGEE & 0.0061 & 0.0284 &0.0233 & 0.0033 & 0.0226& 0.0175\\[3pt]
  &  {$\beta_2$} & PULS & -0.0143 &0.0768 & 0.0401 & 0.0023 & 0.0681 & 0.0417 \\
 & & SGEE & 0.0116 & 0.0275 & 0.0125 & -0.0039 &0.0259 &0.0196\\[3pt]
  &  {$\theta_1$} & PULS & -0.0159 &0.0310 & 0.0252 & 0.0031 & 0.0218 & 0.0168 \\
 & & SGEE & -0.0030 & 0.0083 & 0.0045 & 0.0015 &0.0098 &0.0064\\[3pt]
  &  {$\theta_2$} & PULS & -0.0026 &0.0192 & 0.0112 & 0.0048 & 0.0252 & 0.0200 \\
 & & SGEE & 0.0040 & 0.0200 &0.0133 & 0.0002 & 0.0115&0.0084\\[3pt]
  &  {$\theta_3$} & PULS & 0.0151 &0.0331 & 0.0308 & -0.0067 & 0.0228 & 0.0150 \\
 & & SGEE & 0.0006 & 0.0133 & 0.0083 & -0.0018 &0.0103 &0.0064\\
\hline
\end{tabular*}
\end{table}

To investigate the robustness of the SGEE and PULS estimators to
correlation structure misspecification, we also carry out a simulation
study in which an $\operatorname{AR}(1)$ correlation structure is used in the covariance
matrix estimation detailed in Section~\ref{sec4}, when the true underlying
correlation structure is $\operatorname{ARMA}(1,1)$. Table~\ref{table3} reports the results under
this misspecification. 
The table shows that in the presence of correlation structure
misspecification, SGEE still produces more efficient parameter
estimates than PULS.

We also include a simulated example where the covariates in $\mathbf
{Z}$ follow discrete distributions. The same model as above is used
except that the covariates $\mathbf{X}_i^\top(t_{ij})$ are drawn
independently from a three-dimensional Gaussian distribution with mean
$\mathbf{0}$, variance $1$ and pairwise correlation $0.1$, and
$\mathbf{Z}_i^\top(t_{ij})$ are independently drawn from a binomial
distribution with success probability $0.5$. The errors $e_i(t_{ij})$
are generated with the $\operatorname{AR}(1)$ serial correlation structure of $(\gamma
,\rho)=(1,0.9)$. The simulation results for this example are presented
in Table~\ref{table4}. The same finding as above can be obtained. Some additional
results, that is, those on the average angles between the estimated and
the true parameter vectors, are given in Appendix D of the
supplementary material [\citet{CLLWsupp}].


\subsection{Real data analysis}\label{subsec52}


We next illustrate the partially linear single-index model and the
proposed SGEE estimation method through an empirical example which
explores the relationship between lung function and air pollution.
There is voluminous literature studying the effects of air pollution on
people's health. For a review of the literature, the reader is referred
to \citet{PBR95}. Many studies have found association between air
pollution and health problems such as increased respiratory symptoms,
decreased lung function, increased hospitalizations or hospital visits
for respiratory and cardiovascular diseases and increased respiratory
morbidity [\citet{DSSWSF89}, \citet{KWSDSF89}, \citet{POPE91},
\citet{BASGRW92}, \citet{LH92}]. While
earlier research often used time series or cross-sectional data to
evaluate the health effects of air pollution, recent advances in
longitudinal data analysis techniques offer greater opportunities for
studying this problem. In this paper, we will examine whether air
pollution has a significant adverse effect on lung function, and, if
so, to what extent. The use of the partially linear single-index model
and the SGEE method would provide greater modeling flexibility than
linear models and allow the within-subject correlation to be adequately
taken into account. We will use a longitudinal data set obtained from a
study where a total of 971 4th-grade children aged between 8 and 14
years (at their first visit to the hospital/clinic) were followed over
10 years. For each yearly visit of the children to the hospital/clinic,
records on their forced expiratory volume (FEV), asthma symptom at
visit (\mbox{ASSPM}, 1 for those with symptoms and 0 for those without),
asthmatic status (ASS, 1 for asthma patient and 0 for nonasthma
patient), gender (G, 1 for males and 0 for females), race (R, 1 for
nonwhites and 0 for whites), age (A), height (H), BMI 
and respiratory infection at visit (RINF, 1 for those with infection
and 0 for those without) were taken. Together with the measurements
from the children, the mean levels of ozone and NO$_2$ in the month
prior to the visit were also recorded. Due to dropout or other reasons,
the majority of children had 4 to 5 years of records, and the total
number of observations in the data set is 3809.

As in many other studies, the FEV will be used as a measure of lung
function, and its log-transformed values, log(FEV), will be used as the
response values in our model. Our main interest is to determine whether
higher levels of ozone and NO$_2$ would lead to decrements in lung
function. To account for the effects of other confounding factors, we
include all other recorded variables.
As age and height exhibit strong co-linearity (with a correlation of
0.78), we will only use height in the study. In fitting the partially
linear single-index model to the data, all the continuous variables
(i.e., FEV, H, BMI, OZONE and NO$_2$) are log-transformed, and the
log(BMI), log(OZONE) and log(NO$_2$) are included in the single-index
part. The log(H) and all the binary variables are included in the
linear part of the model.

\begin{figure}

\includegraphics{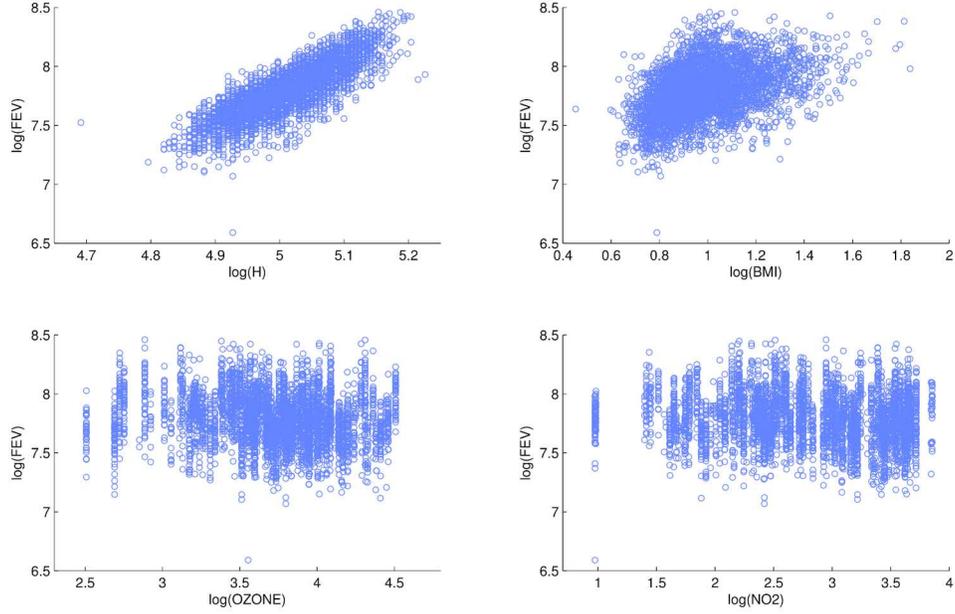}

\caption{The scatter plots of the response variable $\log(\mathrm{FEV})$ against
the continuous regressors, that is, (clockwise from top left) $\log(\mathrm{H})$,
$\log(\mathrm{BMI})$, $\log(\mathrm{NO}_2)$, $\log(\mathrm{OZONE})$.}\vspace*{18pt}\label{fig3}
\end{figure}

\begin{figure}

\includegraphics{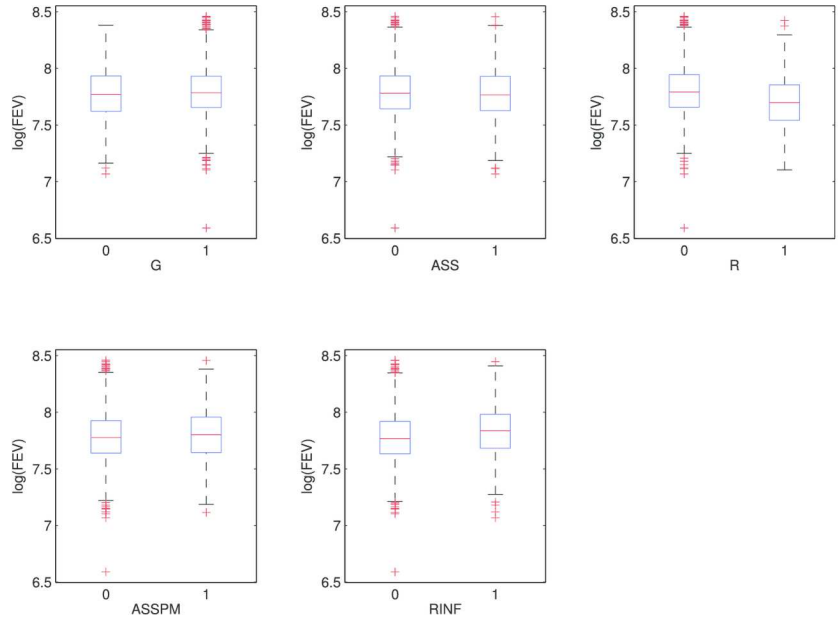}

\caption{The box plots of the response variable $\log(\mathrm{FEV})$ against the
binary regressors, that is, (clockwise from top left) G, ASS, R, RINF, ASSPM.}\label{fig4}
\end{figure}

\begin{figure}[b]

\includegraphics{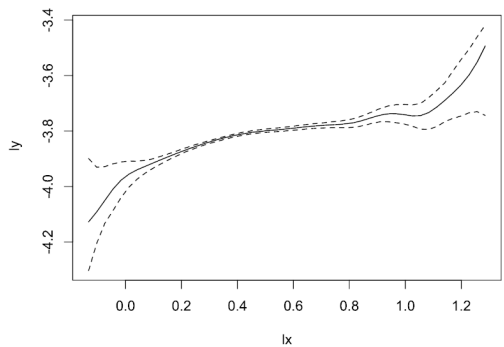}

\caption{The estimated link function and its 95\% point-wise
confidence intervals.}\label{fig5}
\end{figure}

The scatter plots of the response variable against the continuous
regressors are shown in Figure~\ref{fig3}, and the box plots of the response
against the binary regressors are given in Figure~\ref{fig4}. We use an
$\operatorname{ARMA}(1,1)$ within-subject correlation structure in the estimation of the
covariance matrix for the proposed SGEE method. The resulting estimated
model is as follows:
\begin{eqnarray*}
&& \log(\mbox{FEV})
\\
&& \approx 0.0325*\mbox{G}-0.0111*\mbox{ASS}-0.0671*\mbox{R}
\\[-4pt]
&&\hspace*{11pt} (0.0041)\hspace*{22pt} (0.0080)\hspace*{35pt} (0.0059)
\\
&&\quad\qquad{} -0.0047*\mbox{ASSPM}-0.0068*\mbox{RINF}+2.3206*\log(\mbox{H}),
\\[-4pt]
&&\hspace*{33pt}\hspace*{13pt} (0.0085)\hspace*{50pt} (0.0043)\hspace*{40pt} (0.0307)
\\
&&\quad\qquad {} +\widehat{\eta} \bigl[0.9929*\log(\mbox{BMI})-0.0924*\log(\mbox{OZONE})-0.0753*\log(\mbox{NO}_2) \bigr]
\\[-4pt]
&&\hspace*{55pt} (0.0560)\hspace*{58pt} (0.0127)\hspace*{74pt} (0.0125),
\end{eqnarray*}
where the numbers in the parentheses under the estimated coefficien's
are their respective estimated standard errors. The estimated link
function and its 95\% point-wise confidence intervals are plotted in
Figure~\ref{fig5}. 

%

From Figure~\ref{fig5}, it can be seen that the estimated link function is
overall increasing. The $95\%$ point-wise confidence intervals show
that a linear functional form for the unknown link function would be
rejected, and thus the partially liner single-index model might be more
appropriate than the traditional linear regression model. Meanwhile, it
can be seen from the above estimated model that height and BMI are
significant positive factors in accounting for lung function. Taller
children and children with larger BMI tend to have higher FEV.
Furthermore, male and white children have, on average, higher FEV than
female or nonwhite children. Furthermore, both OZONE and NO$_2$ in the
single-index component have negative effects on children's lung
function, as the estimated coefficients for OZONE and NO$_2$ are
negative, and the estimated link function is increasing. Although these
negative effects are relatively small in magnitude compared to the
effect of BMI, they are statistically significant. This means that
higher levels of ozone and NO$_2$ tend to lead to reduced lung function
as represented by lower values of~FEV.

%

\section{Conclusions and discussions}
\label{sec6}

In this paper, we study a partially linear single-index modeling
structure for possibly unbalanced longitudinal data in a general
framework, which includes both the sparse and dense longitudinal data
cases. An SGEE method with the first-stage local linear smoothing is
introduced to estimate the two parameter vectors as well as the
unspecified link function.

In Theorems~\ref{teo1} and~\ref{teo2}, we derive the
asymptotic properties of the proposed parametric and nonparametric
estimators in different scenarios, from which we find that the
convergence rates and asymptotic variances of the resulting estimators
in the sparse longitudinal data case could be substantially different
from those in the dense longitudinal data. In Section~\ref{sec4}, we propose a
semiparametric method to estimate the error covariance matrices which
are involved in the estimation equations. The conditional variance
function is estimated by using the log-transformed local linear method,
and the parameters in the correlation matrices are estimated by the
minimum generalized variance method. In particular, if the correlation
matrices are correctly specified, as is stated in Corollary~\ref{co1}, the
SGEE-based estimators $\widehat{\bolds\beta}$ and $\widehat
{\bolds\theta}$ are generally\vspace*{2pt} asymptotically more efficient than
the corresponding PULS estimators $\widetilde{\bolds\beta}$ and
$\widetilde{\bolds\theta}$ in the sense that the asymptotic
covariance matrix of the SGEE estimators minus that of the PULS
estimators is negative semi-definite. Both the simulation study and
empirical data analysis in Section~\ref{sec5} show that the proposed methods
work well in the finite samples.

Recently, \citet{YL13} developed a new nonparametric regression
function estimation method for a longitudinal regression model. This
method takes into account the within-subject correlation information
and thus generally improves the asymptotic estimation efficiency. It
would also be interesting to incorporate the within-subject correlation
information in the local linear estimation of the unknown link function
in this paper and to examine both theoretical and empirical performance
of the resulting estimator. We will leave this issue for future
research. Another possible future topic is to extend the semiparametric
techniques of variable selection and specification testing proposed by
\citet{LLLT10} from the i.i.d. case to the general longitudinal
data case discussed in the present paper.


\begin{appendix}
\section{Regularity conditions}\label{appA}
To establish the asymptotic properties of the SGEE estimators proposed
in Section~\ref{sec2}, we introduce the following regularity conditions,
although some of
them might not be the weakest possible.

\begin{assumption}\label{assump1}\label{as1}
The kernel function $K(\cdot)$ is a bounded and symmetric probability
density function with compact support. Furthermore, the kernel function
has a continuous first-order derivative function denoted by $\dot
{K}(\cdot)$.
\end{assumption}

\begin{assumption}\label{as2}
(i) The errors ${e}_{ij}\equiv e_i(t_{ij})$, $1\leq i\leq n$, $1\leq
j\leq m_i$, are independent across $i$; that is, ${\mathbf e}_i$
defined in Section~\ref{sec2}, $1\leq i\leq n$, are mutually independent.

(ii) The covariates ${\mathbf X}_{ij}$ and ${\mathbf Z}_{ij}$, $1\leq
i\leq n$, $1\leq j\leq m_i$, are i.i.d. random vectors.


(iii) The errors $e_{ij}$ are independent of the covariates
${\mathbf{Z}}_{ij}$ and ${\mathbf{X}}_{ij}$, and for each~$i$,
$e_{ij}$, $1\leq j\leq m_i$, may be correlated with each other.
Furthermore, $\mathrm{E}[e_{ij}]=0$, $0<\mathrm{E}[e_{ij}^2]<\infty$ and
$\mathrm{E}[\mid e_{ij}\mid^{2+\delta}]<\infty$ for some $\delta>0$. The\vspace*{2pt}
largest eigenvalues of $\mathbf{W}_i$ and $\mathbf{W}_i \mathrm{E}[\mathbf{e}_i\mathbf{e}_i]\mathbf{W}_i$
are bounded for any~$i$.
\end{assumption}

\begin{assumption}\label{as3}
(i) The density function $f_{\bolds\theta}(\cdot)$ of ${\mathbf
X}_{ij}^\top{\bolds\theta}$ is positive and has a continuous
second-order derivative in $\mathcal{U}=\{{\mathbf{x}}^\top
{\bolds\theta}\dvtx  {\mathbf x}\in{\mathcal{X}}, {\bolds\theta
}\in\Theta\}$, where $\Theta$ is a compact parameter space for
${\bolds\theta}$ and ${\mathcal{X}}$ is a compact support of
${\mathbf{X}}_{ij}$.

(ii) The function $\rho_{\mathbf{Z}}(u\mid{\bolds\theta})=\mathrm{E}[{\mathbf Z}_{ij}\mid{\mathbf X}_{ij}^\top{\bolds\theta}=u]$ has a
bounded and continuous second-order derivative (with respect to $u$)
for any ${\bolds\theta}\in\Theta$, and $\mathrm{E} [\llVert
{\mathbf{Z}}_{ij}\rrVert^{2+\delta} ]<\infty$, where $\delta$ was
defined in {Assumption~\ref{as2}(iii)}.
\end{assumption}

\begin{assumption}\label{as4}
The link function $\eta(\cdot)$ has continuous derivatives up to the
second order.
\end{assumption}

\begin{assumption}\label{as5}
The bandwidth $h$ satisfies
%
\begin{equation}
\label{eqA1} \omega_nh^6\rightarrow0,\qquad
\frac{n^2h^2}{N_{n}(h)\log n}\rightarrow\infty,\qquad \frac{T_n^{\afrac
{2}{2+\delta}}\log n}{h^2N_n(h)}=o(1),
\end{equation}
where\vspace*{1pt} $N_n(h)=\sum_{i=1}^n{1}/{(m_ih)}$, $T_n=\sum_{i=1}^nm_i$ and
$\delta$ was defined in Assumption~\ref{as2}(iii). Furthermore,
$\max_{1\leq i \leq n} (m_i^4+m_i^3h^{-1})=o(w_n)$.
\end{assumption}

\begin{remark}\label{re5}
Assumption~\ref{as1} imposes some mild restrictions on the kernel functions,
which have been used in the existing literature in i.i.d. and weakly
dependent time series cases; see, for example, \citet{FG96}
and \citet{Gao07}. The compact support restriction on the kernel functions
can be removed if we impose certain restrictions on the tail of the
kernel function. In Assumption~\ref{as2}(i), the longitudinal data under
investigation is assumed to be independent across subjects $i$, which
is not uncommon in longitudinal data analysis; see, for example, \citet{WZ06}
and \citet{ZFS}. Assumption~\ref{as2}(ii) is imposed to
simplify the presentation of the asymptotic results. However, we may
replace Assumption~\ref{as2}(ii) with the conditions that the covariates
${\mathbf X}_{ij}$ and ${\mathbf Z}_{ij}$ are i.i.d. across $i$ and
identically distributed across $j$, and in the case of dense
longitudinal data, it is further satisfied that for $\kappa
=0,1,2,\ldots,$
%
\begin{equation}
\label{eqA11} \operatorname{Var} \Biggl[\frac{1}{m_i}\sum
_{j=1}^{m_i} \frac{U_{ij}}{h} \biggl(
\frac{{\mathbf X}_{ij}^\top{\bolds\theta}-u}{h} \biggr)^{\kappa} K
\biggl(\frac{{\mathbf X}_{ij}^\top{\bolds\theta}-u}{h} \biggr)
\Biggr]\leq C(m_ih)^{-1}
\end{equation}
uniformly for $u\in{\mathcal U}$ and ${\bolds\theta}\in\Theta
$, where $U_{ij}$ can be $1$, ${\mathbf Z}_{ij}B_1({\mathbf Z}_{ij})$,
or ${\mathbf X}_{ij}B_2({\mathbf X}_{ij})$, $B_1(\cdot)$ and
$B_2(\cdot)$ are two bounded functions, and $C$ is a positive constant
which is independent of $i$. When ${\mathbf X}_{ij}$ and ${\mathbf
Z}_{ij}$ are stationary and $\alpha$-mixing dependent across $j$ for
the case of dense longitudinal data, it is easy to validate the
high-level condition (\ref{eqA11}). In Assumption~\ref{as2}(iii), we allow
the error terms to have certain within-subject correlation, which makes
the model assumptions more realistic.\vadjust{\goodbreak} Assumption~\ref{as3} gives some
commonly-used conditions in partially linear single-index models; see
\citet{XH06} and \citet{CGL13b}, for example.
Assumption~\ref{as4} is a mild smoothness condition on the link function
imposed for the application of the local linear fitting. Assumption~\ref{as5}
gives a set of restrictions on the bandwidth $h$, which is involved in
the estimation of the link function. Note that the bandwidth conditions
in Assumption~\ref{as5} imply that the milder bandwidth conditions in (C.1) of
Lemma 1 in the supplemental material [\citet{CLLWsupp}] are
satisfied. Hence we can use Lemma 1 to prove our main theoretical results.
\end{remark}

We next give some regularity conditions, which are needed to derive the
asymptotic property of the nonparametric conditional variance
estimators in Section~\ref{sec4}.

\begin{assumption}\label{as6}
The kernel function $K_1(\cdot)$ is a continuous and symmetric
probability density function with compact support.
\end{assumption}

\begin{assumption}\label{as7}
The observation times, $t_{ij}$, are i.i.d. and have a continuous and
positive probability density function $f_T(t)$, which has a compact
support~${\mathcal{T}}$. The density function of $\xi^2(t_{ij})$ is
continuous and bounded. Let $\delta>2$, which strengthens the moment
conditions in Assumptions~\ref{as2} and~\ref{as3}.
\end{assumption}

\begin{assumption}\label{as8}
The conditional variance function $\sigma^2(\cdot)$ has a continuous
second-order derivative and satisfies $\inf_{t\in{\mathcal
{T}}}\sigma^2(t)>0$. Let $\dot{\sigma}^2(\cdot)$ and $\ddot{\sigma
}^2(\cdot)$ be its first-order and second-order derivative functions,
respectively.
\end{assumption}

\begin{assumption}\label{as9}
The bandwidth $h_1$ satisfies
%
\begin{equation}
\label{eqA2} h_1\rightarrow0,\qquad \frac{T_n^{\afrac{2}{2+\delta/2}}\log
n}{h_1^2N_n(h_1)}=o(1),
\end{equation}
where $N_n(h_1)=\sum_{i=1}^n{1}/{(m_ih_1)}$.
\end{assumption}

\begin{remark}\label{re6}
Assumption~\ref{as7} imposes a mild condition on the observation times [see,
e.g., \citet{JW11}] and strengthens the moment conditions on
$e_{ij}$ and ${\mathbf Z}_{ij}$. However, such moment conditions are
not uncommon in the asymptotic theory for nonparametric conditional
variance estimation [\citet{CCP09}]. Since the local linear
smoothing technique is applied, a certain smoothness condition has to
be assumed on $\sigma^2(\cdot)$, as is done in Assumption~\ref{as8}.
Assumption~\ref{as9} gives some mild restrictions on the bandwidth $h_1$, which
is used in the estimation of the conditional variance function.
\end{remark}


\section{Proofs of the main results}\label{appB}

In this appendix, we provide the detailed proofs of the main results
given in Section~\ref{sec3}.

\subsection{Proof of Theorem \texorpdfstring{\protect\ref{teo1}}{1}}

By the definition of the weighted local linear estimators in (\ref
{eq24}) and (\ref{eq25}), we have
%
\begin{eqnarray}\label{eqB1}
\widehat{\eta}(u\mid{\bolds\beta}, {\bolds\theta})-\eta
(u)&=&\sum
_{i=1}^n{\mathbf s}_i(u\mid{
\bolds\theta}) ({\mathbf Y}_i-{\mathbf{Z}}_i{
\bolds\beta} )-\eta(u)
\nonumber
\\
&=&\sum_{i=1}^n{\mathbf
s}_i(u\mid{\bolds\theta}){\mathbf e}_i+\sum
_{i=1}^n{\mathbf s}_i(u\mid{
\bolds\theta}){\mathbf{Z}}_i ({\bolds\beta}_0-{
\bolds\beta} )
\nonumber
\\
&&{}+\sum_{i=1}^n{\mathbf
s}_i(u\mid{\bolds\theta}) \bigl[{\bolds\eta}({
\mathbf{X}}_i,{\bolds\theta}_0)-{\bolds\eta}({\mathbf{X}}_i,{\bolds\theta}) \bigr]
\\
&&{}+\sum_{i=1}^n{\mathbf
s}_i(u\mid{\bolds\theta}){\bolds\eta}({
\mathbf{X}}_i,{\bolds\theta})-\eta(u)
\nonumber
\\
&\equiv& I_{n1}+I_{n2}+I_{n3}+I_{n4}.\nonumber
\end{eqnarray}

For $I_{n1}$, note that by a first-order Taylor expansion of $K(\cdot
)$, we have, for $i=1,\dots,n$ and $j=1,\dots,m_i$,
\[
K \biggl(\frac{{\mathbf X}_{ij}^\top{\bolds\theta}-u}{h} \biggr)=K
\biggl(\frac{{\mathbf X}_{ij}^\top{\bolds\theta}_0-u}{h} \biggr)+ \dot{K}
\biggl(\frac{{\mathbf
X}_{ij}^\top{\bolds\theta}_*-u}{h} \biggr)\frac{{\mathbf
{X}}_{ij}^\top({\bolds\theta}-{\bolds\theta}_0)}{h},
\]
where $\dot{K}(\cdot)$ is the first-order derivative of $K(\cdot)$
and
${\bolds\theta}_*={\bolds\theta}_0+\lambda
_*({\bolds\theta}-{\bolds\theta}_0)$,
$0<\lambda_*<1$. Hence, by some standard calculations and the
assumption that ${n^2h^2}/{\{N_n(h)\log n\}}\rightarrow\infty$, we
have
%
\begin{eqnarray}\label{eqB2}
I_{n1}&=&\sum_{i=1}^n{\mathbf
s}_i(u\mid{\bolds\theta}_0){\mathbf
e}_i+\sum_{i=1}^n \bigl[{
\mathbf s}_i(u\mid{\bolds\theta})-{\mathbf s}_i(u
\mid{\bolds\theta}_0) \bigr]{\mathbf e}_i
\nonumber
\\
&=&\sum_{i=1}^n{\mathbf
s}_i(u\mid{\bolds\theta}_0){\mathbf
e}_i+O_P \biggl(\llVert{\bolds\theta}-{
\bolds\theta}_0\rrVert\cdot\frac{\sqrt{N_n(h)\log n}}{nh} \biggr)
\\
&=&\sum_{i=1}^n{\mathbf
s}_i(u\mid{\bolds\theta}_0){\mathbf
e}_i+o_P \bigl(\llVert{\bolds\theta}-{\bolds\theta}_0\rrVert\bigr)\nonumber
\end{eqnarray}
for any $u\in{\mathcal{U}}$ and ${\bolds\theta}\in\Theta$.

By Lemma 2 in the supplementary material [\citet{CLLWsupp}], we can
prove that
%
\begin{equation}
\label{eqB3} I_{n2}=-\rho_{\mathbf{Z}}^\top(u) ({
\bolds\beta}-{\bolds\beta}_0)+O_P \bigl(
\llVert{\bolds\beta}-{\bolds\beta}_0\rrVert
^2+\llVert{\bolds\theta}-{\bolds\theta}_0
\rrVert^2 \bigr)
\end{equation}
for any $u\in{\mathcal{U}}$, where $\rho_{\mathbf{Z}}(u)\equiv\rho
_{\mathbf{Z}}(u\mid{\bolds\theta}_0)=\mathrm{E}[{\mathbf
{Z}}_{ij}\mid{\mathbf{X}}_{ij}^{\top}{\bolds\theta}_0=u]$.

Note that
\[
\eta\bigl({\mathbf{X}}_{ij}^\top{\bolds\theta}\bigr)-
\eta\bigl({\mathbf{X}}_{ij}^\top{\bolds\theta}_0\bigr)=\dot{\eta}\bigl({\mathbf{X}}_{ij}^\top{
\bolds\theta}_0\bigr) {\mathbf{X}}_{ij}^\top({
\bolds\theta}-{\bolds\theta}_0)+O_P\bigl(
\llVert{\bolds\theta}-{\bolds\theta}_0\rrVert
^2\bigr),
\]
which, together with Lemma 3 in the supplementary material [\citet{CLLWsupp}], leads to
%
\begin{equation}
\label{eqB4} I_{n3}=-\dot{\eta}(u)\rho_{\mathbf{X}}^\top(u)
({\bolds\theta}-{\bolds\theta}_0)+O_P \bigl(
\llVert{\bolds\theta}-{\bolds\theta}_0\rrVert
^2 \bigr)
\end{equation}
for any $u\in{\mathcal{U}}$, where $\rho_{\mathbf{X}}(u)\equiv\rho
_{\mathbf{X}}(u\mid{\bolds\theta}_0)=\mathrm{E} [{\mathbf
{X}}_{ij}\mid{\mathbf{X}}_{ij}^{\top}{\bolds\theta}_0=u ]$.\vspace*{1pt}

By a second-order Taylor expansion of $\eta(\cdot)$ and the
first-order Taylor expansion of $K(\cdot)$ used to handle $I_{n1}$, we
can prove that, for any $u\in{\mathcal{U}}$, we have
%
\begin{equation}
\label{eqB5} I_{n4}=\tfrac{1}{2}\mu_2\ddot{
\eta}(u)h^2\bigl[1+O_P(h)\bigr]+o_P\bigl(
\llVert{\bolds\theta}-{\bolds\theta}_0\rrVert\bigr).
\end{equation}

Recall that $\widehat{\bolds\beta}$ and $\widehat{\bolds\theta}_1$ are the solutions to the equations in (\ref{eq27}). By
(\ref{eqB1})--(\ref{eqB5}), we can prove that, uniformly for
$i=1,\dots,n$ and $j=1,\dots,m_i$,
\begin{eqnarray}\label{eqB6}
&& \widehat{\eta} \bigl({\mathbf{X}}_{ij}^\top\widehat{
\bolds\theta}_1\mid\widehat{\bolds\beta}, \widehat{
\bolds\theta}_1 \bigr)-\eta\bigl({\mathbf{X}}_{ij}^\top
{\bolds\theta}_0\bigr)
\nonumber
\\[-2pt]
&&\qquad =\widehat{\eta} \bigl({\mathbf{X}}_{ij}^\top\widehat{
\bolds\theta}_1\mid\widehat{\bolds\beta}, \widehat{
\bolds\theta}_1 \bigr)-\widehat{\eta} \bigl({\mathbf
{X}}_{ij}^\top{\bolds\theta}_0\mid
\widehat{\bolds\beta}, \widehat{\bolds\theta}_1 \bigr)+
\widehat{\eta}\bigl({\mathbf{X}}_{ij}^\top{\bolds\theta}_0\mid\widehat{\bolds\beta}, \widehat{\bolds\theta}_1\bigr) -\eta\bigl({\mathbf{X}}_{ij}^\top{
\bolds\theta}_0\bigr)
\nonumber
\\[-2pt]
&&\qquad =\widehat{\dot{\eta}} \bigl({\mathbf{X}}_{ij}^\top{
\bolds\theta}_0\mid\widehat{\bolds\beta}, \widehat{
\bolds\theta}_1 \bigr){\mathbf{X}}_{ij}^\top
(\widehat{\bolds\theta}_1-{\bolds\theta}_0 )+
\widehat{\eta}\bigl({\mathbf{X}}_{ij}^\top{\bolds\theta}_0 \mid\widehat{\bolds\beta}, \widehat{\bolds\theta}_1\bigr) -\eta\bigl({\mathbf{X}}_{ij}^\top{
\bolds\theta}_0\bigr)\nonumber
\\[-2pt]
&&\quad\qquad{} +O_P \bigl(\llVert\widehat{\bolds\theta}_1-{
\bolds\theta}_0\rrVert^2 \bigr)
\\[-2pt]
&&\qquad=\dot{\eta}\bigl({\mathbf{X}}_{ij}^\top{\bolds\theta}_0\bigr) \bigl[{\mathbf{X}}_{ij}-\rho_{\mathbf{X}}
\bigl({\mathbf{X}}_{ij}^\top{\bolds\theta}_0
\bigr) \bigr]^\top(\widehat{\bolds\theta}_1-{
\bolds\theta}_0 ) \bigl(1+o_P(1)\bigr)\nonumber
\\[-2pt]
&&\quad\qquad{} +\sum
_{k=1}^n{\mathbf s}_k\bigl({
\mathbf{X}}_{ij}^\top{\bolds\theta}_0\bigr){
\mathbf e}_k
-\rho_{\mathbf{Z}}^\top\bigl({\mathbf{X}}_{ij}^\top{
\bolds\theta}_0\bigr) (\widehat{\bolds\beta}-{\bolds\beta}_0) \bigl(1+o_P(1)\bigr)\nonumber
\\[-2pt]
&&\quad\qquad{} + \tfrac{1}{2}
\mu_2\ddot{\eta}\bigl({\mathbf{X}}_{ij}^\top{
\bolds\theta}_0\bigr)h^2+O_P
\bigl(h^3\bigr)
+O_P \bigl(\llVert\widehat{\bolds\theta}_1-{
\bolds\theta}_0\rrVert^2+\llVert\widehat{
\bolds\beta}-{\bolds\beta}_0\rrVert^2
\bigr),\nonumber
\end{eqnarray}
where\vspace*{1pt} ${\mathbf s}_k({\mathbf{X}}_{ij}^\top{\bolds\theta
}_0)\equiv{\mathbf s}_k({\mathbf{X}}_{ij}^\top{\bolds\theta
}_0\mid{\bolds\theta}_0)$.\vspace*{1pt}

By the definitions of $\widehat{\bolds\beta}$ and $\widehat
{\bolds\theta}_1$ [see (\ref{eq27}) in Section~\ref{sec2}], we have
%
\begin{equation}
\label{eqB7} \sum_{i=1}^n\widehat{
\bolds\Lambda}_i^\top(\widehat{\bolds\theta}_1){\mathbf{W}}_i \bigl[{\mathbf
Y}_i-{\mathbf{Z}}_i\widehat{\bolds\beta}-
\widehat{\bolds\eta}({\mathbf{X}}_i\mid\widehat{\bolds\beta}, \widehat{\bolds\theta}_1) \bigr]={\mathbf{0}}.
\end{equation}
By the uniform consistency results for the local linear estimators
(such as Lemmas 2 and 3 in the supplementary material [\citet{CLLWsupp}]), we can approximate $\widehat{\bolds\Lambda
}_i(\widehat{\bolds\theta}_1)$ in (\ref{eqB7}) by
${\bolds\Lambda}_i={\bolds\Lambda}_i({\bolds\theta
}_0)$ when deriving the asymptotic distribution theory. Then we have
%
\begin{eqnarray}\label{eqB8}
{\mathbf0}&=&\sum_{i=1}^n\widehat{\bolds\Lambda}_i^\top(\widehat{\bolds\theta}_1){
\mathbf{W}}_i \bigl[{\mathbf Y}_i-{\mathbf
{Z}}_i\widehat{\bolds\beta}- \widehat{\bolds\eta}({
\mathbf{X}}_i\mid\widehat{\bolds\beta},\widehat{\bolds\theta}_1) \bigr]
\nonumber
\\[-2pt]
&=&\sum_{i=1}^n{\bolds\Lambda}_i^\top{\mathbf{W}}_i \bigl[{\mathbf
Y}_i-{\mathbf{Z}}_i\widehat{\bolds\beta}-
\widehat{\bolds\eta}({\mathbf{X}}_i\mid\widehat{\bolds\beta},\widehat{\bolds\theta}_1) \bigr]
\\[-2pt]
&&{}+\sum_{i=1}^n\bigl(\widehat{\bolds\Lambda}_i(\widehat{\bolds\theta}_1)-{\bolds\Lambda}_i\bigr)^\top{\mathbf{W}}_i \bigl[{
\mathbf Y}_i-{\mathbf{Z}}_i\widehat{\bolds\beta}-
\widehat{\bolds\eta}({\mathbf{X}}_i\mid\widehat{\bolds\beta},\widehat{\bolds\theta}_1) \bigr]
\nonumber
\\[-2pt]
&\stackrel{P} {\sim}& \sum_{i=1}^n{
\bolds\Lambda}_i^\top{\mathbf{W}}_i
\bigl[{\mathbf Y}_i-{\mathbf{Z}}_i\widehat{\bolds\beta}- \widehat{\bolds\eta}({\mathbf{X}}_i\mid\widehat{
\bolds\beta},\widehat{\bolds\theta}_1) \bigr]
\bigl[1+O_P\bigl(\llVert\widehat{\bolds\theta}_1-{
\bolds\theta}_0\rrVert\bigr) \bigr],\nonumber
\end{eqnarray}
where and below $a_n\stackrel{P}\sim b_n$ denotes $a_n=b_n(1+o_P(1))$.
Furthermore, note that
\begin{eqnarray}
{\mathbf Y}_i-{\mathbf{Z}}_i\widehat{\bolds\beta}- \widehat{\bolds\eta}({\mathbf{X}}_i\mid\widehat{
\bolds\beta},\widehat{\bolds\theta}_1)={\mathbf
e}_i-{\mathbf{Z}}_i (\widehat{\bolds\beta} -{
\bolds\beta}_0 )- \bigl[\widehat{\bolds\eta}({
\mathbf{X}}_i\mid\widehat{\bolds\beta},\widehat{\bolds\theta}_1)- {\bolds\eta}({\mathbf{X}}_i,{
\bolds\theta}_0) \bigr],
\nonumber
\end{eqnarray}
which, together with (\ref{eqB6}) and the bandwidth condition $\omega
_nh^6=o(1)$, implies that
\begin{eqnarray}\label{eqB9}
&& \sum_{i=1}^n{\bolds\Lambda}_i^\top{\mathbf{W}}_i \bigl[{\mathbf
Y}_i-{\mathbf{Z}}_i\widehat{\bolds\beta}-
\widehat{\bolds\eta}({\mathbf{X}}_i\mid\widehat{\bolds\beta},\widehat{\bolds\theta}_1) \bigr]
\nonumber
\\
&&\qquad =\sum_{i=1}^n{\bolds\Lambda}_i^\top{\mathbf{W}}_i{\mathbf
e}_i-\sum_{i=1}^n{\bolds\Lambda}_i^\top{\mathbf{W}}_i{\mathbf
{Z}}_i (\widehat{\bolds\beta} -{\bolds\beta}_0 )\nonumber
\\
&&\quad\qquad{} -\sum_{i=1}^n{
\bolds\Lambda}_i^\top{\mathbf{W}}_i
\bigl[\widehat{\bolds\eta}({\mathbf{X}}_i\mid\widehat{
\bolds\beta},\widehat{\bolds\theta}_1)- {\bolds\eta}({
\mathbf{X}}_i,{\bolds\theta}_0) \bigr]
\nonumber
\\
&&\qquad =-\sum_{i=1}^n{\bolds\Lambda}_i^\top{\mathbf{W}}_i \bigl[{
\mathbf{Z}}_i-{\bolds\rho}_{\mathbf{Z}}({\mathbf
{X}}_{i},{\bolds\theta}_0) \bigr] (\widehat{
\bolds\beta} -{\bolds\beta}_0 ) \bigl(1+o_P(1)
\bigr)
\nonumber\\[-8pt]\\[-8pt]\nonumber
&&\quad\qquad{} -\sum_{i=1}^n{\bolds\Lambda}_i^\top{\mathbf{W}}_i \bigl\{ \bigl[
\dot{\bolds\eta}({\mathbf{X}}_i,{\bolds\theta
}_0)\otimes{\mathbf1}_{p}^\top\bigr]\odot
\bigl[{\mathbf{X}}_i- {\bolds\rho}_{\mathbf{X}}({
\mathbf{X}}_{i},{\bolds\theta}_0) \bigr] \bigr\}\nonumber
\\
&&\hspace*{58pt}{} \times (
\widehat{\bolds\theta}_1 -{\bolds\theta}_0 )
\bigl(1+o_P(1)\bigr)
\nonumber
\\
&&\quad\qquad{} +\sum_{i=1}^n{\bolds\Lambda}_i^\top{\mathbf{W}}_i \Biggl[{\mathbf
e}_i-\sum_{k=1}^n{\mathbf
s}_{k}({\mathbf{X}}_i,{\bolds\theta}_0){
\mathbf e}_k \Biggr]\nonumber
\\
&&\quad\qquad{} +O_P \bigl(\llVert\widehat{
\bolds\beta}-{\bolds\beta}_0\rrVert^2+\llVert
\widehat{\bolds\theta}_1-{\bolds\theta}_0
\rrVert^2 \bigr),\nonumber
\end{eqnarray}
where
${\mathbf s}_{k}({\mathbf{X}}_i,{\bolds\theta}_0)=
[{\mathbf s}_{k}^\top({\mathbf{X}}_{i1}^\top{\bolds\theta
}_0),\dots,
{\mathbf s}_{k}^\top({\mathbf{X}}_{im_i}^\top{\bolds\theta
}_0) ]^\top$, ${\bolds\rho}_{\mathbf{Z}}({\mathbf
{X}}_{i},{\bolds\theta}_0)$ and ${\bolds\rho}_{\mathbf
{X}}({\mathbf{X}}_{i},{\bolds\theta}_0)$ were defined in
Section~\ref{sec2}. Following the standard proof in the existing literature
[see, e.g., \citet{Ich93}, \citet{CGL13b}], we can show the weak
consistency of $\widehat{\bolds\beta}$ and $\widehat
{\bolds\theta}_1$. Note that
\begin{eqnarray}
&& \sum_{i=1}^n{\bolds\Lambda}_i^\top{\mathbf{W}}_i{\bolds\Lambda}_i \pmatrix{ \widehat{\bolds\beta}-{\bolds\beta}_0
\vspace*{3pt}\cr
\widehat{\bolds\theta}_1-{\bolds\theta}_0}
\nonumber
\\
&&\qquad =\sum_{i=1}^n{\bolds\Lambda}_i^\top{\mathbf{W}}_i \bigl\{ \bigl[
\dot{\bolds\eta}({\mathbf{X}}_i,{\bolds\theta
}_0)\otimes{\mathbf1}_{p}^\top\bigr]\odot
\bigl[{\mathbf{X}}_i- {\bolds\rho}_{\mathbf
X}({
\mathbf{X}}_{i},{\bolds\theta}_0) \bigr] \bigr\} (
\widehat{\bolds\theta}_1 -{\bolds\theta}_0 )
\nonumber
\\
&&\quad\qquad{} +\sum_{i=1}^n{\bolds\Lambda}_i^\top{\mathbf{W}}_i \bigl[{
\mathbf{Z}}_i-{\bolds\rho}_{\mathbf
Z}({\mathbf{X}}_{i},{
\bolds\theta}_0) \bigr] (\widehat{\bolds\beta} -{
\bolds\beta}_0 )
\nonumber
\end{eqnarray}
and
\[
\sum_{i=1}^n{\bolds\Lambda}_i^\top{\mathbf{W}}_i \Biggl[\sum
_{k=1}^n{\mathbf s}_{k}({
\mathbf{X}}_i,{\bolds\theta}_0){\mathbf
e}_k \Biggr]=o_P\bigl(\omega_n^{1/2}
\bigr), %
\]
which, together with (\ref{eqB8}) and (\ref{eqB9}), lead to
%
\begin{eqnarray}
\Biggl[\sum_{i=1}^n{\bolds\Lambda}_i^\top{\mathbf{W}}_i{\bolds\Lambda}_i \Biggr]
\pmatrix{ \widehat{\bolds\beta}-{\bolds\beta}_0
\vspace*{3pt}\cr
\widehat{\bolds\theta}_1-{\bolds\theta}_0}\stackrel{P}\sim\sum_{i=1}^n{
\bolds\Lambda}_i^\top{\mathbf{W}}_i{
\mathbf e}_i.\label{eqB10}
\end{eqnarray}

Define ${\mathbf I}({\bolds\theta}_0, {\mathbf B}_0)=\operatorname{diag} \{{\mathbf I}_d, {\mathbf M} \}, {\mathbf
O}({\bolds\theta}_0)={
{\mathbf O}_{d\times d}\ \ {\mathbf O}_{d\times1}\choose
{\mathbf O}_{p\times d}\ \ {\bolds\theta}_0}$,
where ${\mathbf M}=({\bolds\theta}_0, {\mathbf B}_0)$ was
defined in Section~\ref{sec3}. It is easy to find that
%
\begin{equation}
\label{eqB11} {\mathbf I}_{d+p}={\mathbf I}({\bolds\theta}_0, {\mathbf B}_0){\mathbf I}^\top({
\bolds\theta}_0, {\mathbf B}_0)={\mathbf O}({
\bolds\theta}_0){\mathbf O}^\top({\bolds\theta
}_0)+{\mathbf I}({\mathbf B}_0){\mathbf
I}^\top({\mathbf B}_0).
\end{equation}
By the identification condition on ${\bolds\theta}_0$, we may
show that
\begin{eqnarray}
\widehat{\bolds\theta}-{\bolds\theta}_0&=&
\frac
{\widehat{\bolds\theta}_1}{\llVert \widehat{\bolds\theta
}_1\rrVert }-\frac{{\bolds\theta}_0}{\llVert {\bolds\theta
}_0 \rrVert }=\frac{\widehat{\bolds\theta}_1}{\llVert \widehat
{\bolds\theta}_1\rrVert }-\frac{{\bolds\theta}_0}{\llVert
\widehat{\bolds\theta}_1\rrVert }+
\frac{{\bolds\theta}_0}{\llVert \widehat{\bolds\theta
}_1\rrVert }-\frac{{\bolds\theta}_0}{\llVert {\bolds\theta
}_0 \rrVert }
\nonumber
\\
&\stackrel{P}\sim&\frac{\widehat{\bolds\theta
}_1-{\bolds\theta}_0}{\llVert {\bolds\theta}_0\rrVert }-
{\bolds\theta}_0{
\bolds\theta}_0^\top\frac{\widehat
{\bolds\theta}_1-{\bolds\theta}_0}{\llVert {\bolds\theta}_0\rrVert }= \bigl({\mathbf
I}_p-{\bolds\theta}_0{\bolds\theta}_0^\top\bigr) (\widehat{\bolds\theta}_1-{\bolds\theta}_0 ),
\nonumber
\end{eqnarray}
which implies that $\widehat{\bolds\theta}-{\bolds\theta
}_0={\mathbf B}_0{\mathbf B}_0^\top(\widehat{\bolds\theta
}_1-{\bolds\theta}_0 )$
and
%
\begin{equation}
\label{eqB12} \pmatrix{ \widehat{\bolds\beta}-{\bolds\beta}_0
\vspace*{3pt}\cr
\widehat{\bolds\theta}-{\bolds\theta}_0}={\mathbf I}({\mathbf B}_0){\mathbf
I}^\top({\mathbf B}_0)\pmatrix{ \widehat{\bolds\beta}-{\bolds\beta}_0
\vspace*{3pt}\cr
\widehat{\bolds\theta}_1-{\bolds\theta}_0}.
\end{equation}

By (\ref{eqB10}), (\ref{eqB11}) and using the fact that
${\bolds\Lambda}_i{\mathbf O}({\bolds\theta}_0)={\mathbf0}$,
we have
\[
{\mathbf I}^\top({\mathbf B}_0) \Biggl[\sum
_{i=1}^n{\bolds\Lambda}_i^\top{
\mathbf{W}}_i{\bolds\Lambda}_i \Biggr]{\mathbf I}({
\mathbf B}_0){\mathbf I}^\top({\mathbf B}_0)
\pmatrix{ \widehat{\bolds\beta}-{\bolds\beta}_0
\vspace*{3pt}\cr
\widehat{\bolds\theta}_1-{\bolds\theta}_0}
\stackrel{P}\sim{\mathbf I}^\top({\mathbf B}_0) \Biggl[
\sum_{i=1}^n{\bolds\Lambda
}_i^\top{\mathbf{W}}_i{\mathbf
e}_i \Biggr], %
\]
which, together with (\ref{eqB12}), implies that
\begin{eqnarray}
\pmatrix{ \widehat{\bolds\beta}-{\bolds\beta}_0
\vspace*{3pt}\cr
\widehat{\bolds\theta}-{\bolds\theta}_0} &\stackrel{P}
\sim&{\mathbf I}({\mathbf B}_0) \Biggl\{ {\mathbf I}^\top({
\mathbf B}_0) \Biggl[\sum_{i=1}^n{
\bolds\Lambda}_i^\top{\mathbf{W}}_i{
\bolds\Lambda}_i \Biggr]{\mathbf I}({\mathbf B}_0)
\Biggr\}^{-1}{\mathbf I}^\top({\mathbf B}_0)
\Biggl[\sum_{i=1}^n{\bolds\Lambda}_i^\top{\mathbf{W}}_i{\mathbf
e}_i \Biggr].
\nonumber
\end{eqnarray}
Thus, by (\ref{eq31})--(\ref{eq321}), the definition of the
Moore--Penrose inverse and the classical central limit theorem for
independent sequence, we can show that (\ref{eq33}) in Theorem~\ref{teo1}
holds.

\subsection{Proof of Corollary \texorpdfstring{\protect\ref{co1}}{1}}

By Theorem~\ref{teo1}, the PULS estimators $\widetilde{\bolds\beta}$ and
$\widetilde{\bolds\theta}$ have the following asymptotic normal
distribution:
%
\begin{equation}
\label{eqB13} \omega_n^{1/2}\pmatrix{ \widetilde{
\bolds\beta}-{\bolds\beta}_0
\vspace*{3pt}\cr
\widetilde{\bolds\theta}-{\bolds\theta}_0} \stackrel{d}\longrightarrow\mathrm{N}
\bigl({\mathbf0}, {\bolds\Omega}_{0\ast}^{+}{\bolds\Omega}_{1\ast}{\bolds\Omega}_{0\ast}^{+} \bigr),
\end{equation}
where ${\bolds\Omega}_{0\ast}$ and ${\bolds\Omega
}_{1\ast}$ are two matrices such that
\[
\frac{1}{\omega_n}\sum_{i=1}^n{\bolds\Lambda}_i^\top{\bolds\Lambda}_i
\stackrel{P}\rightarrow{\bolds\Omega}_{0\ast},\qquad \frac{1}{\omega
_n}\sum
_{i=1}^n\mathrm{E} \bigl[{\bolds\Lambda}_i^\top{\mathbf{V}}_i{\bolds\Lambda}_i \bigr]\rightarrow{\bolds\Omega}_{1\ast},
\]
and ${\mathbf{V}}_i$ is the conditional covariance matrix of ${\mathbf
{e}}_i$.

On the other hand, when the weights ${\mathbf{W}}_i$, $i=1,\dots,n$,
are chosen as the inverse of ${\mathbf{V}}_i$, by Theorem~\ref{teo1}, we have
%
\begin{equation}
\label{eqB14} \omega_n^{1/2} \pmatrix{ \widehat{\bolds\beta}-{\bolds\beta}_0
\vspace*{3pt}\cr
\widehat{\bolds\theta}-{
\bolds\theta}_0} \stackrel{d}\longrightarrow\mathrm{N} \bigl({
\mathbf0}, {\bolds\Omega}_{\ast}^{+} \bigr),
\end{equation}
where
${\bolds\Omega}_{\ast}$ is a positive semi-definite matrix such that
\[
\frac{1}{\omega_n}\sum_{i=1}^n\mathrm{E}
\bigl[{\bolds\Lambda}_i^\top{\mathbf{V}}_i^{-1}{
\bolds\Lambda}_i \bigr]\rightarrow{\bolds\Omega}_{\ast}.
\]

In order to prove Corollary~\ref{co1}, by (\ref{eqB13}) and (\ref{eqB14}),
we need only to show
${\bolds\Omega}_{0\ast}^{+}{\bolds\Omega}_{1\ast
}{\bolds\Omega}_{0\ast}^{+}-{\bolds\Omega}_\ast^+$ is
positive semi-definite. Letting ${\bolds\Theta}_i={\bolds\Omega}_{0\ast}^+{\bolds\Lambda}_i{\mathbf
{V}}_i^{1/2}-{\bolds\Omega}_\ast^+{\bolds\Lambda
}_i{\mathbf{V}}_i^{-1/2}$, we have
\begin{eqnarray}
{\bolds\Theta}_i{\bolds\Theta}_i^\top&=&
\bigl({\bolds\Omega}_{0\ast}^+{\bolds\Lambda}_i{
\mathbf{V}}_i^{1/2}-{\bolds\Omega}_\ast^+{
\bolds\Lambda}_i{\mathbf{V}}_i^{-1/2}
\bigr) \bigl({\bolds\Omega}_{0\ast
}^+{\bolds\Lambda}_i{
\mathbf{V}}_i^{1/2}-{\bolds\Omega}_\ast^+{
\bolds\Lambda}_i{\mathbf{V}}_i^{-1/2}
\bigr)^\top
\nonumber
\\
&=&{\bolds\Omega}_{0\ast}^+{\bolds\Lambda}_i{
\mathbf{V}}_i{\bolds\Lambda}_i{\bolds\Omega}_{0\ast
}^+-{\bolds\Omega}_{0\ast}^+{\bolds\Lambda
}_i{\bolds\Lambda}_i{\bolds\Omega}_{\ast}^+-{\bolds\Omega}_{\ast}^+{\bolds\Lambda}_i{\bolds\Lambda}_i{\bolds\Omega}_{0\ast}^++{\bolds\Omega}_{\ast
}^+{\bolds\Lambda}_i{\mathbf{V}}_i^{-1}{\bolds\Lambda}_i{\bolds\Omega}_{\ast}^+,
\nonumber
\end{eqnarray}
which indicates that
%
\begin{equation}
\label{eqB15} \frac{1}{\omega_n}\sum_{i=1}^n
\mathrm{E} \bigl[{\bolds\Theta}_i{\bolds\Theta}_i^\top\bigr]\rightarrow{\bolds\Omega
}_{0\ast}^{+}{\bolds\Omega}_{1\ast}{\bolds\Omega}_{0\ast}^{+}-{\bolds\Omega}_\ast^+.
\end{equation}
As $\mathrm{E} [{\bolds\Theta}_i{\bolds\Theta}_i^\top
]$ is positive semi-definite, by (\ref{eqB15}) we know that
${\bolds\Omega}_{0\ast}^{+}{\bolds\Omega}_{1\ast
}{\bolds\Omega}_{0\ast}^{+}-{\bolds\Omega}_\ast^+$ is
also positive semi-definite. Hence the proof of Corollary~\ref{co1} is complete.

\subsection{Proof of Theorem \texorpdfstring{\protect\ref{teo2}}{2}}

Note that
\begin{eqnarray}\label{eqB16}
\widehat{\eta}(u)-\eta(u)&=&\sum_{i=1}^n{
\mathbf s}_i(u \mid\widehat{\bolds\theta}) \bigl({\mathbf
Y}_i-{\mathbf{Z}}_i^\top\widehat{\bolds\beta} \bigr)-\eta(u)
\nonumber
\\
&=&\sum_{i=1}^n{\mathbf
s}_i(u\mid\widehat{\bolds\theta}){\mathbf e}_i+
\Biggl[\sum_{i=1}^n{\mathbf
s}_i(u\mid\widehat{\bolds\theta}){\bolds\eta}({
\mathbf{X}}_i,{\bolds\theta}_0)-\eta(u) \Biggr]
\nonumber\\[-8pt]\\[-8pt]\nonumber
&&{}+\sum_{i=1}^n{\mathbf
s}_i(u\mid\widehat{\bolds\theta}) {\mathbf{Z}}_i^\top
({\bolds\beta}_0-\widehat{\bolds\beta} )
\nonumber
\\
&\equiv& I_{n1,\ast}+I_{n2,\ast}+I_{n3,\ast}.\nonumber
\end{eqnarray}

By Assumption~\ref{as1}, we have
%
\begin{equation}
\label{eqB17} \qquad K \biggl(\frac{{\mathbf X}_{ij}^\top\widehat{\bolds\theta
}-u}{h} \biggr)=K \biggl(\frac{{\mathbf X}_{ij}^\top{\bolds\theta
}_0-u}{h}
\biggr)+ \dot{K} \biggl(\frac{{\mathbf X}_{ij}^\top{\bolds\theta
}_\lozenge
-u}{h} \biggr)\frac{{\mathbf{X}}_{ij}^\top(\widehat{\bolds\theta}-{\bolds\theta}_0)}{h},
\end{equation}
where ${\bolds\theta}_\lozenge={\bolds\theta}_0+\lambda
_\lozenge(\widehat{\bolds\theta}-{\bolds\theta}_0)$ for
some $0<\lambda_\lozenge<1$. By Theorem~\ref{teo1}, we have
%
\begin{equation}
\label{eqB18} \llVert\widehat{\bolds\theta}-{\bolds\theta}_0\rrVert+\llVert\widehat{\bolds\beta}-{\bolds\beta}_0\rrVert=O_P\bigl(\omega_n^{-1/2}
\bigr).
\end{equation}
It follows from (\ref{eqB17}), (\ref{eqB18}) and (\ref{eq340}) that
%
\begin{eqnarray}\label{eqB19}
\qquad I_{n3,\ast}&=&\sum_{i=1}^n{\mathbf
s}_i(u\mid{\bolds\theta}_0) {
\mathbf{Z}}_i^\top({\bolds\beta}_0-
\widehat{\bolds\beta} )+\sum_{i=1}^n
\bigl[{\mathbf s}_i(u\mid\widehat{\bolds\theta})-{\mathbf
s}_i(u\mid{\bolds\theta}_0) \bigr]{
\mathbf{Z}}_i^\top({\bolds\beta}_0-
\widehat{\bolds\beta} )\nonumber
\\
&=&O_P\bigl(\omega_n^{-1/2}
\bigr)+O_P\bigl(\omega_n^{-1}
\bigr)
\\
&=& o_P\bigl(\varphi_n^{-1/2}(h)
\bigr).\nonumber
\end{eqnarray}

Similar to the proof of (\ref{eqB5}), we can show that
%
\begin{equation}
\label{eqB20} I_{n2,\ast}=\tfrac{1}{2}\ddot{\eta}(u)
\mu_2h^2\bigl(1+o_P(1)\bigr).
\end{equation}

For $I_{n1,\ast}$, note that by (\ref{eqB17}) and (\ref{eqB18}),
we can show that $\sum_{i=1}^n{\mathbf s}_i(u\mid{\bolds\theta
}_0){\mathbf e}_i$ is the leading term of $I_{n1,\ast}$. Letting
$z_i({\bolds\theta}_0)={\mathbf s}_i(u\mid{\bolds\theta
}_0){\mathbf e}_i$ and by Assumption~\ref{as2}, it is easy to check that $\{
z_i({\bolds\theta}_0)\dvtx  i\geq1\}$ is a sequence of independent
random variables. By Assumption~\ref{as2}(iii), we have $\mathrm{E}[z_i({\bolds\theta}_0)]=0$. By (\ref{eq340}), (\ref
{eq34}) and the central limit theorem, it can be readily seen that\vspace*{-2pt}
%
\begin{equation}
\label{eqB21} \varphi_n^{1/2}(h) I_{n1,\ast}
\stackrel{d}\rightarrow\mathrm{N}\bigl(0,\sigma_*^2\bigr).
\end{equation}
In view of (\ref{eqB16}), (\ref{eqB19})--(\ref{eqB21}), the proof
of Theorem~\ref{teo2} is complete. 
\end{appendix}





\section*{Acknowledgements}
The authors wish to thank the
Co-editor, the Associate Editor and two referees for their valuable
comments and suggestions, which substantially improved an earlier
version of the paper.

\begin{supplement}[id=suppA]
\stitle{Supplement to ``Semiparametric GEE analysis in partially linear single-index models for longitudinal data''}
\slink[doi]{10.1214/15-AOS1320SUPP} 
\sdatatype{.pdf}
\sfilename{aos1320\_supp.pdf}
\sdescription{The supplement gives the proof of Theorem~\ref{teo3} and some technical lemmas
that were used to prove the main results in Appendix \ref{appB}. It also
includes some additional results of our simulation studies described in
Section~\ref{sec5}.}
\end{supplement}


\printaddresses
\end{document}